\documentclass[11pt,a4paper]{article}

\usepackage{amssymb,latexsym,amsmath,graphics}

\usepackage[left=2.7cm, right=2.7cm, top=2.7cm]{geometry}

\usepackage{float}
\usepackage{xcolor}
\usepackage[demo]{graphicx}
\usepackage[boxed]{algorithm2e}
\usepackage{graphicx,color}
\usepackage{latexsym}
\usepackage{epsfig}
\usepackage{amssymb}
\usepackage{amstext}
\usepackage{amsgen}
\usepackage{amsxtra}
\usepackage{amsgen}
\usepackage{amsthm}
\usepackage{tocloft} 
\usepackage{authblk}
\usepackage{graphicx}
\usepackage{caption}
\usepackage{subcaption}
\usepackage{enumerate}  
\usepackage{enumitem}
\usepackage{color}
\usepackage{dsfont}
\usepackage{amsmath}
\usepackage{geometry}
\usepackage{marginnote}
\usepackage{bbm}
\usepackage{cite}
\usepackage{authblk}
\usepackage{hyperref} 
\usepackage{mathtools}
\usepackage{authblk}

\newlist{abbrv}{itemize}{1}
\setlist[abbrv,1]{label=,labelwidth=1in,align=parleft,itemsep=0.1\baselineskip,leftmargin=!}

\numberwithin{equation}{section}

\newcommand{{\Dd}}{\mathcal{D}}
\newcommand{\Tt}{\Gamma}

\def\vk{k}

\def \tt {{\widetilde t}}

\def\tilde{\widetilde}
\def \Noise {\mathbb{W}}

\def \prob{{\mathbb P}}
\def \Pp{{\mathbb P}}
\def \bra{\langle}
\def \cet{\rangle}

\def \a {\beta}

\def \e {\varepsilon}

\newcommand{\argmin}{\mathop{\mathrm{arg\,min}}}
\newcommand{\Nn}{{\mathcal N}}
\newcommand{\E}{{\mathbb E}}
\newcommand{\F}{{\mathcal F}}
\newcommand{\N}{{\mathbb N}}
\newcommand{\R}{{\mathbb R}}
\newcommand{\Z}{{\mathbb Z}}

\newcommand{\expec}{{\mathbb E \,}}

\newcommand{\g}{{\, |\,}}
\newcommand{\oT}{\textbf{T}}

\newtheorem {theorem}{Theorem}
\newtheorem {proposition}[theorem]{Proposition}
\newtheorem {definition}[theorem]{Definition}
\newtheorem {lemma}[theorem]{Lemma}
\newtheorem {corollary}[theorem]{Corollary}

\begin {document}

\title{
Random tree Besov priors -- Towards fractal imaging}

\author[1]{Hanne Kekkonen\thanks{Email: h.n.kekkonen@tudelft.nl}}
\author[2]{Matti Lassas}
\author[2]{Eero Saksman}
\author[2]{Samuli Siltanen}
\affil[1]{Delft Institute of Applied Mathematics, TU Delft, The Netherlands}
\affil[2]{Department of Mathematics and Statistics, University of Helsinki, Finland}
\date{\today}                    
\setcounter{Maxaffil}{0}

\maketitle

\textbf{Abstract.} We propose alternatives to Bayesian a priori distributions that are frequently used in the study of inverse problems. Our aim is to construct priors that have similar good edge-preserving properties as total variation or Mumford-Shah priors but correspond to well defined infinite-dimensional random variables, and can be approximated by finite-dimensional random variables. We introduce a new wavelet-based model, where the non zero coefficient are chosen in a systematic way so that prior draws
have certain fractal behaviour. We show that realisations of this new prior take values in some Besov spaces and have  singularities only on a small set $\tau$ that has a certain Hausdorff dimension. We also introduce an efficient algorithm for calculating the MAP estimator, arising from the the new prior, in denoising problem.

\section{Introduction}

\noindent
Inverse problems arise from the need to extract information from indirect measurements. They are typically ill-posed, meaning that algorithmic recovery of information is sensitive to noise and modelling errors. Robust reconstruction methods are based on combining the measurement data with {\it a priori} knowledge about the unknown target.  
Formulating {\it a priori} knowledge mathematically 
is a core challenge in inverse problems research. Popular 
models for {\it a priori} information  promote global smoothness, piecewise regularity, or sparsity in a given or learned basis or a more general collection of building blocks. They can be implemented using variational regularisation, providing a stable solution but little information about the uncertainties of the model. Bayesian inversion offers an attractive alternative, 
by also delivering information on how uncertainties in the data or model affect obtained point estimates. 

Practical measurements are always finite and corrupted by noise, which can in many cases be reasonably modelled using independent Gaussian random variables. This gives rise to a discrete measurement model of the form
\begin{equation}\label{measmodel}
M_i=(Af)_i+ w_i, \quad i=1,\dots,n, \ w_i\overset{\textnormal{iid}}{\sim} \Nn(0,1),
\end{equation}
where $A$ describes the forward process. Computational solution of the inverse problem requires also a finite approximate model for the unknown $f$. 
It is advisable to design the models for $f$ and for the noise in a discretisation-invariant way \cite{Haario2004a, Lasanen2007, Lassas2004, Lassas2009}. One way to achieve this is to build a continuous model which is discretised as late as possible in the analysis and solution procedure \cite{Dashti2016, Stuart2010}. In this paper we consider the case where $A:H_1\to H_2$ is a linear operator between Banach spaces $H_1$ and $H_2$, and assume the continuous equivalent model (in the sense of \cite{Brown1996, Reiss2008})
\begin{align}\label{eq:IP}
M=Af+\varepsilon\Noise,
\end{align}
where $\Noise$ is a Gaussian white noise process indexed by $H_2$.

Fractals have many applications in science since they often describe the real world better than traditional mathematical models. As Mandelbrot noted, "Clouds are not spheres, mountains are not cones, coastlines are not circles, and bark is not smooth, nor does lightning travel in a straight line. Fractal geometry and analysis of fractal dimension is a powerful tool that has been used to study turbulent mixing flows and atmospheric turbulence \cite{Davis1996, Ribak1994}, evaluate the risk of a patient developing cancer and predict tumour malignancy  \cite{Li2007, Chan2016, Caldwell1990}, screening for osteopenia \cite{Alman2012}, recognising Alzheimer's disease patients from magnetoencephalogram recordings \cite{Gomez2009}, and in study of macular diseases \cite{Thomas2014}.

\begin{figure}[h]
\begin{picture}(400,120)
\put(0,110){(a)}
\put(20,0){\includegraphics[width=6.5cm]{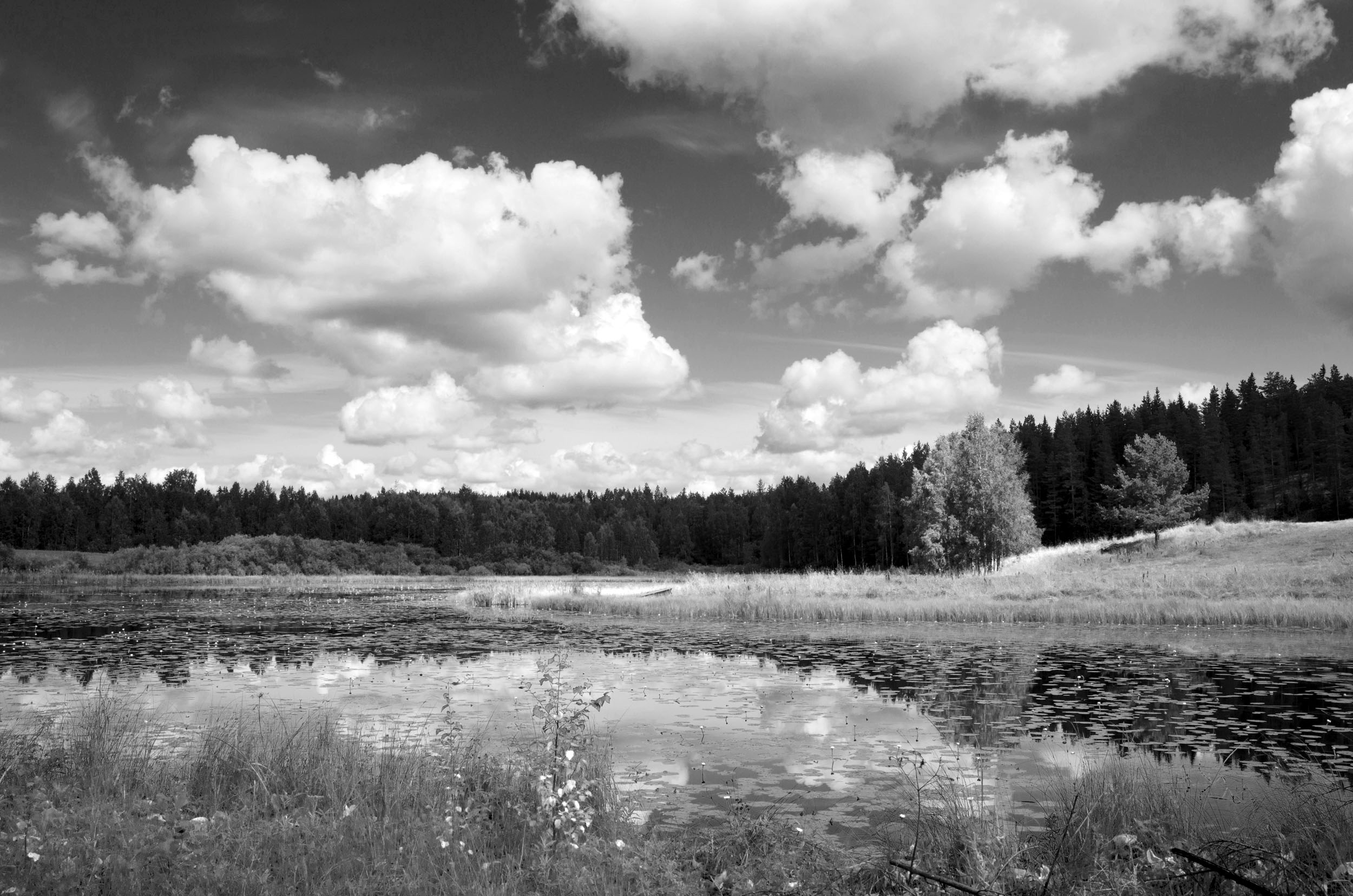}}
\put(300,0){\includegraphics[width=5cm]{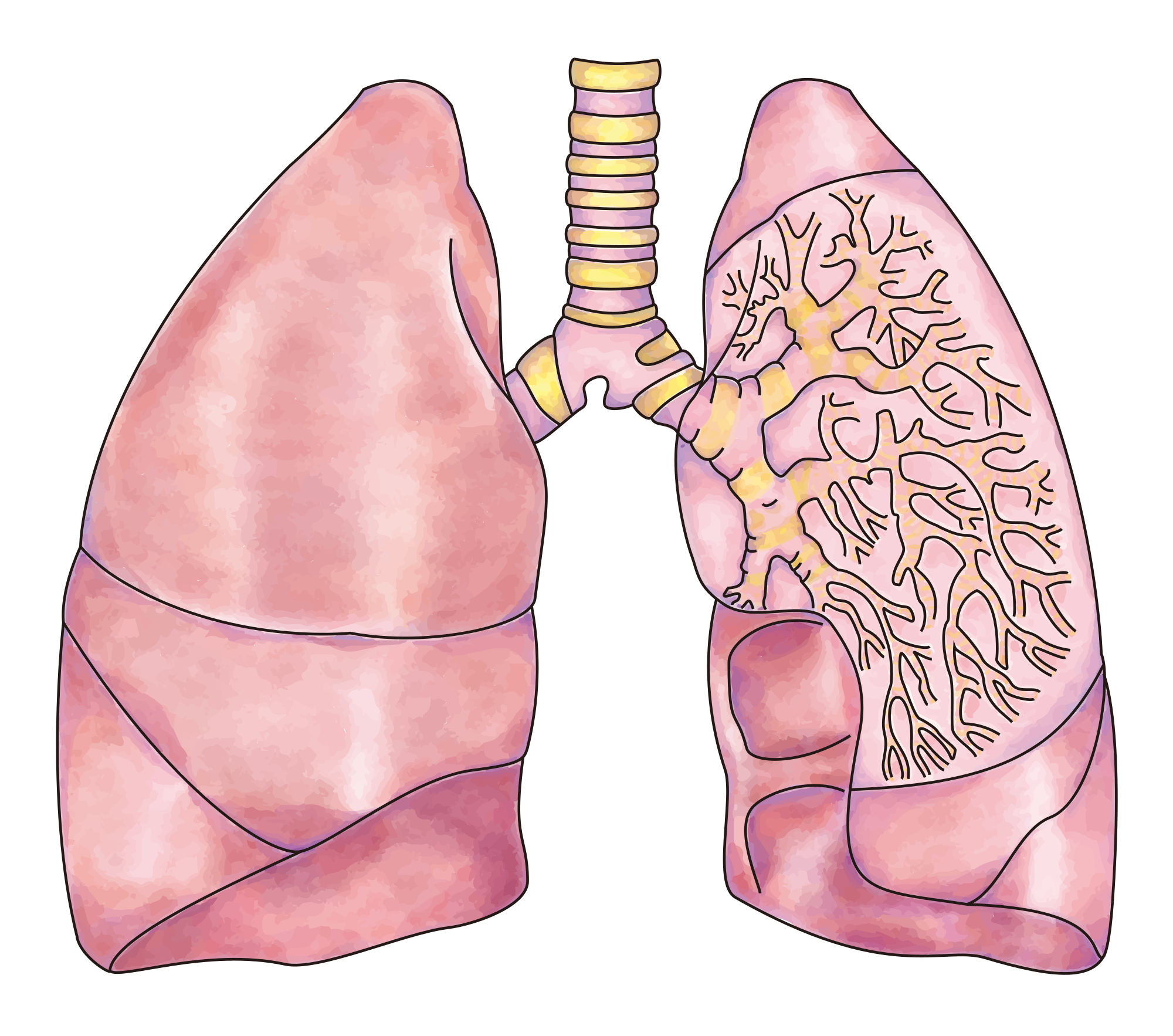}}
\put(240,36){\includegraphics[width=3.5cm]{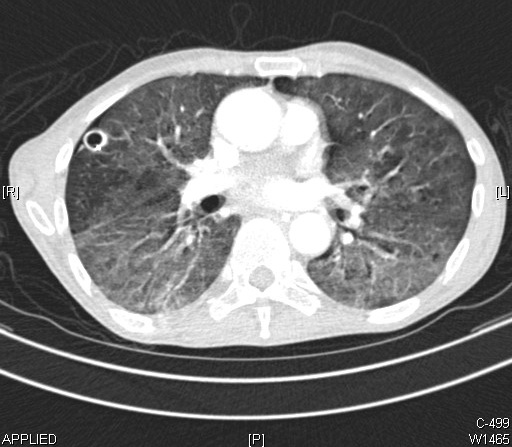}}
\put(220,110){(b)}
\end{picture}
\caption{\label{fig:fractalexamples}(a) Photograph featuring fractal-like structures such as boundaries of clouds. (b) The three-dimensional structure of human lungs follows a self-similar rule leading to a fractal dimension. Also shown is an X-ray tomography slice of lungs. Images courtesy of Wikimedia Commons.}
\end{figure}

We will next give a few examples of measurement models like (\ref{eq:IP}) arising in practical applications where the unknown quantity has fractal properties.
In each case we encounter a function $f$ that has singular support on a fractal set $\tau$.
 We emphasise that in the paper, we mostly consider random fractals, that are sets produced by a random process which (Hausdorff) dimension is not an integer. 
As a one dimensional example we consider functions $f:[0,1]\rightarrow\R$ whose graphs have fractal properties or which smoothness changes on a fractal set. Such functions appear, for example, in medicine (heart rate, stride, and breathing variability) \cite{west2012fractal} and financial data (local volatility models for asset prices where the volatility changes when the price drops below a threshold value) \cite{mandelbrot2005inescapable,mandelbrot2005parallel,LejayPigato2019,pigato2019}. Here the formula (\ref{eq:IP}) can model e.g. interpolation or extrapolation problems. 
An example of two-dimensional problems of the form (\ref{eq:IP}) are digital image processing tasks such as deblurring, inpainting and denoising an image $f:[0,1]^2\rightarrow \R$. If the underlying picture is for instance a photograph containing clouds, coastlines or skylines with forest or mountains, the boundaries of the objects are fractals.
 See Figure \ref{fig:fractalexamples} (a).
X-ray tomography is a good example of a three-dimensional inverse problem that fits to our model. There one records several two-dimensional X-ray images of a patient along different directions of projection. The goal is to reconstruct the X-ray attenuation coefficient $f:[0,1]^3\rightarrow \R$ from these images, where each pixel is considered as a line integral of $f$ along a ray. The structure of human lungs is known to be fractal  \cite{west2012fractal}, so our new model is well suited for lung tomography. See Figure \ref{fig:fractalexamples} (b).

Most of the theory for infinite-dimensional Bayesian inverse problems focuses on the case were the unknown is assumed to follow a Gaussian prior. Gaussian inverse problems benefit from fast computational properties but in many signal and image reconstruction problems the detection of edges and interfaces is crucial and these problems are poorly modelled by Gaussian priors. One way of circumventing the issue is to use hierarchical Gaussian models as in  \cite{Calvetti2007, Calvetti2008 }. Closely related to this paper are the hierarchical models whose maximum a posteriori estimate converges to a minimiser of the Mumford-Shah functional \cite{Helin2009, Helin2011}. We follow the idea of using wavelet-based Besov priors introduced in \cite{Lassas2009} and further studied in \cite{Dashti2012}. Consistency of such priors has also been considered in \cite{Agapiou2018}.
These priors are especially useful since smooth functions with few local irregularities have sparser expansion in the wavelet basis than, e.g., in the Fourier basis. If the unknown is assumed to be sparse in some basis, using a prior that encourages this sparsity results in more efficient finite-dimensional approximation of the solution. Methods for recovering finite-dimensional estimates of the unknown based on wavelet bases are broadly studied in image and signal processing, and statistical literature, see e.g. \cite{Abramovich1998, Chambolle1998, Daubechies2004, Donoho1994}.

We adopt the Bayesian approach to inverse problems and assign a prior probability measure $\Pi$ to $f$. The solution to the Bayesian  inverse problem is the conditional distribution of $f$ given data, and the mean or mode of the posterior can be used as a point estimator. One often used method for achieving edge-preserving solutions in image analysis is to employ the so-called total variation prior and use the mode of the posterior as a point estimate. In practice this means solving the minimisation problem
\begin{align}\label{minimisation2}
\min_{f\in H} \big\{\|Af-M\|_{L^2}^2+\beta \|\nabla f\|_{L^1}\big\},
\end{align}  
where $H$ is a finite-dimensional subspace of piece-wise smooth functions  and $A:L^2\to L^2$ is a continuous operator. 
Despite of active research on this area, no natural infinite
dimensional models have been found that would 
result (\ref{minimisation2}) as a maximum a posteriori estimate.
In other words, the widely used formal prior 
\begin{align}\label{eq:formalPrior}
\pi(f)\underset{formally}\propto\exp(-\a \|\nabla f\|_{L^1})
\end{align}
is not known to correspond  any well defined random variable. 
It is also known that the usual discrete total variation priors can converge to Gaussian smoothness priors when discretisation becomes denser, see  \cite{Lassas2004}.   
 
The idea presented in \cite{Lassas2009} is to replace formula \eqref{eq:formalPrior} by a well defined prior
\begin{align*}
\pi(f)\underset{formally}\propto\exp(- \|\nabla f\|_{B^0_{pp}}^p),
\end{align*}
where the Besov spaces $B^0_{pp}$ are closely related to $L^p$ spaces. Hence the Besov $B^0_{11}$-priors have similar properties to total variation prior but correspond to well defined infinite-dimensional random variables. The construction of the Besov prior is done through a generalisation of the Karhunen-Lo\`eve expansion and can hence be approximated with finite-dimensional random variables.

In this paper we introduce a new wavelet-based computational model for  {\it a priori} information about the fractal dimension of the unknown target. 
We will build on the idea from \cite{Lassas2009} but choose the non-zero wavelet coefficients in the Karhunen-Lo\`eve expansion in a systematic way, so that the resulting priors have a certain fractal dimension. This is done by introducing a new  random variable $T$ that takes values in the space of `trees' and using a prior
\begin{align*}
\pi(f,T)\underset{formally}\propto\exp(-(C_T^{-1} f,f)_{L^2})\, \pi_{T}(T) 
\end{align*}
where $C_T$ is a covariance matrix depending on the hyper-variable $T$. We choose the tree $T$ so that the set where the Schwartz kernel of $C_T$ is non-smooth has almost surely a low Hausdorff dimension, which means that the realisations $f$ have singularities only on a small set. This opens up new possibilities in medical imaging and signal processing. Our new model allows rigorous analysis in the Bayesian inversion framework. Also, we present an efficient algorithm for the computation of maximum a posteriori (MAP) estimates for two- and three-dimensional denoising problems. 

The rest of the paper is organised as follows. In Section \ref{Sec:Besov} we introducing the Besov space setting that is used for constructing the priors.  We then define the new random tree Besov priors and formulate the main results of the paper in Section \ref{Sec:Fractals}. Section \ref{Sec:exaples} is dedicated for the denoising examples and constructing an algorithm for calculating MAP estimators.

\section{Priors in Besov-spaces}\label{Sec:Besov}

\subsection{Gaussian inverse problems}

We start by motivating the use of semi-Gaussian Besov priors by first considering the standard Gaussian prior. Let $H_1$ is a Banach space, and let $f$ be a $H_1$-valued random variable following a Gaussian distribution $\mu$. We denote by $\Noise$ the Gaussian white noise proces. For simplicity, we assume that $\E f=0$ and denote by $C_f$  and $C_w$ the covariance operators of $f$ and $\Noise$, respectively. If the forward operator $A:H_1\to H_2$ is assumed to be linear the posterior distribution $\mu(\cdot\g M)$ is also Gaussian. It then follows that the conditional mean estimate for the Bayesian inverse problem  \eqref{eq:IP} coincides with the maximum a posteriori estimate (under mild assumptions on $A$, see e.g. \cite{Dashti2013}) and is given by 
\begin{align}\label{ATAformula}
f_{MAP}=(A^*C_w^{-1}A+C_f^{-1})^{-1}A^*C_w^{-1}M, 
\end{align}
where $A^*:H_2\to H_1$ is the adjoint of $A$ (see e.g. \cite{Lehtinen1989}). 


{\bf Example 1.} Here we consider a stereotype of 1-dimensional
linear Bayesian inverse problem, where we use a priori model
that $f$ is a Brownian bridge on interval $(0,\pi)$.
Assume that $H_1=L^2(0,\pi)$ with the usual inner 
product $\bra\cdotp,\cdot\cet$ and basis
 $e_j(x)=\frac 12\sin (jx)$. 
Let $A_N=\{f\in H_1:\ a_j<\bra f,e_j\cet<b_j, \ j=1,2,\dots,N\}$ be a 
cylinder and define $\mu$ by setting 
\begin{align}\label{mu-definition}
\mu(f: (\bra f,e_j\cet)_{j=1}^N\in A_N\})
=c_N\int_{a_1}^{b_1}\dots \int_{a_N}^{b_N} \prod_{j=1}^N
e^{-\lambda_j^{-1}t_j^2/2}dt_1\dots dt_N,
\end{align}
where $c_N$ are normalisation constants and 
$\lambda_j=j^{-2}$. 
We see that $f$ has a representation 
$f=\sum_{j=1}^\infty \lambda_j^{-1/2}f_je_j$ where $f_j$
are independent normalised Gaussian variables. 
Furthermore $f$ has the covariance operator $C_f$ defined by
\begin{align}\label{definition of covariance}
\bra C_f\phi,\psi\cet=\expec(\bra \phi,f\cet\bra f,\psi\cet),
\end{align}
where $\phi,\psi\in L^2(0,\pi)$,
which can easily be seen to be 
\begin{align}\label{C-oper}
C_f\phi=\sum_{j=1}^\infty \lambda_j\bra \phi,e_j\cet e_j.
\end{align}
Now, if we consider the Dirichlet Laplacian, 
that is, the  operator $-\Delta=-\frac {d^2}{dx^2}$ defined in
domain ${\Dd}(-\Delta)=H^2(0,\pi)\cap H^1_0(0,\pi)$
we see that $e_j$ are the complete system of  orthonormal eigenvectors
of $-\Delta $ corresponding to eigenvalues $\lambda_j^{-1}$. 
Thus we see from (\ref{C-oper}) that the covariance operator $C$
coincides with $-\Delta^{-1}$.
In applications distribution (\ref{mu-definition}) is (non-rigorously) expressed by saying that
$f$ has the probability density function
\begin{align}\label{gaussian 1}
\pi_f(f)= ce^{-\bra f,\Delta f\cet/2},\quad f\in L^2(0,\pi)
\end{align}
where $c$ is a normalisation constant. This notation
becomes rigorous if we discretise
the system, that is, consider $N+1$-dimensional vector $(f(0),f(\pi/N),$ $\dots,f(\pi))$ and approximate $-\Delta$ with the finite difference
operator \cite{Lasanen2002}. When $A:L^2(0,\pi)\to L^2(0,1)$ is a bounded linear operator and $\Noise$ is normalised Gaussian 
white noise on interval $(0,1)$  the conditional mean estimate for inverse problem (\ref{eq:IP}), given a measurement $m$, is 
$f_{MAP}=(A^*A-\frac 12\Delta )^{-1}A^*m,$ in accordance with (\ref{ATAformula}).

\subsection{Sobolev and Besov spaces}  

Our aim is to combine the fast computational properties of the Gaussian inverse problems with the good edge-preserving properties of some non-Gaussian variables, in particular those generated  by assuming a total variation prior. Analogously to (\ref{eq:formalPrior}) we would like to use prior
\begin{align*}
\pi(f)\underset{formally}\propto\exp(- \|\nabla f\|_{L^p}^p), \quad p\geq1.
\end{align*}
As this is not a well defined object for general $p$, we will replace it by 
\begin{align*}
\pi(f)\underset{formally}\propto\exp(- \|\nabla f\|_{B^0_{pp}}^p)
\end{align*}
or more generally by 
\begin{align*}
\pi(f)\underset{formally}\propto\exp(- \|f\|_{B^s_{pp}}^p),
\end{align*}
where the Besov space $B^0_{pp}$ has many similar properties
to $L^p$ as we see below. 

In this section we recall the  Besov spaces $B^s_{pq}(\R^d)$, $1\leq p,q<\infty$ and $s\in\R$, and the closely related Sobolev spaces $W^{s,p}(\R^d)$ (see e.g. \cite{Triebel2008} for general theory). 
Let $\F$ denote the Fourier transform and $\mathcal{S'}(\R^d)$ the space of tempered generalised functions (also called tempered distributions). The Sobolev spaces on $\R^d$  for $p\in (1,\infty)$ and $s\in\R$ are defined as 
\begin{align*}
W^{s,p}(\R^d)=\{f\in \mathcal{S'}(\R^d)\, :\, \mathcal{F}^{-1}\big((1+|\cdot|^2)^{\frac{s}{2}}\mathcal{F}f\big)\in L^p(\R^d)\}
\end{align*}
with norm 
\begin{align*}
\|f\|_{W^{s, p}(\mathbb{R}^{d})} :=\big\|\mathcal{F}^{-1}\big((1+|\cdot|^2)^{\frac{s}{2}}\mathcal{F}f\big)\big\|_{L^{p}(\mathbb{R}^{d})}.
\end{align*}
For  $s\in\N$ the above definition is equivalent to having derivatives of $f$ of order at most $s$, defined in the distributional sense,  in $L^p$.  This second definition can be used to consider also the Lebegue exponent $p=1$. For example, $W^{1,1}(\R^d)$ consist of functions that together with their first derivatives lie in $L^1(\R^n)$.

The Besov spaces are slightly more complicated. One possible way
to define them is to use Fourier multipliers, that is,
frequency band filters in the Fourier space.
First, let $\phi_0\in C_0^\infty(\R^d)$ be a function for which
\begin{align*}
\operatorname{supp} (\phi_0) \subset\left\{\xi \in \R^d :|\xi| \leq3/2\right\}, \quad \phi_0(\xi)=1 \text { if }|\xi| \leq 1.
\end{align*}
We can then define 
\begin{align*}
\phi_{j}(\xi)=\phi_0\left(2^{-j} \xi\right)-\phi_0\left(2^{-j+1} \xi\right), \quad \xi \in \R^d,\,\, j\geq 1, 
\end{align*}
and note that 
\begin{align*}
\operatorname{supp} (\phi_{j}) \subset\left\{\xi \in \R^d : 2^{j-1} \leq|\xi| \leq 2^{j+1}\right\}\quad\text{and} 
\quad \sum_{j=0}^{\infty} \phi_{j}(\xi)=1 \text { if } \xi \in \R^d,
\end{align*}
that is, $\phi_j$ is a partition of unity. 
We define the Fourier multipliers $F_j:\mathcal{S}'(\R^d)\to \mathcal{S}'(\R^d)$ as,
\begin{align*}  
F_jf=\F^{-1}(\phi_j(\cdot)\F f).
\end{align*}
The Besov space is the collection of all $f\in\mathcal{S}'$ such that 
\begin{align}\label{eq:BesovNorm1}
\|f\|_{B^s_{pq}(\R^d)}=
\left(\sum_{j=0}^\infty 2^{sjq}\|F_jf\|_{L^p(\R^d)}^q\right)^{1/q} <\infty.
\end{align}

In particular we are interested in $B^1_{1\,1}(\R^d)$ having norm
\begin{align}\label{eq:B111}
\| f\|_{B^1_{1\,1(\R^d)}}
\approx  \|f\|_{B^0_{1\,1}(\R^d)} +\|\nabla f\|_{B^0_{1\,1}(\R^d)}  
\approx \|f\|_{B^0_{1\,1}(\R^d)} +
\sum_{j=0}^\infty \|\nabla(F_jf)\|_{L^1(\R^d)},
\end{align}
as this space is relatively close to space of functions with
bounded variations, $\|\nabla f\|_{L^1}<\infty$. For example, it is not difficult to show that (locally) all the functions $x^\alpha_+$ belong to $ B^1_{1\,1}(\R)$ for any $\alpha>0$.
Thus $B^1_{1\,1}(\R)$  almost contains functions having
Heaviside-type jump singularities. 

We finally mention that Besov norm for $p>1$ can also be defined  equivalently via suitable integrals of finite differences of the function $f$, see e.g. \cite[p.8]{Triebel1992}.

\subsection{Besov spaces valued random variables} 
In this paper we will use the wavelet representation of Besov norms in dimension $d\geq 1$.  For this, we recall that for any integer $r\geq 1$ there exists
compactly supported $C^r$ functions $\phi$ and $\Psi^\ell$  ($\ell=1,\ldots ,2^d-1$) which generate wavelets suitable for multi-resolution analysis of smoothness $C^r$.  More specifically, if we denote

\newcommand{\kauno}{{\mathcal L}}
\begin{align*}
\psi^\ell_{jk}(x)&=2^{jd/2}\Psi^\ell(2^jx-\vk),\qquad \vk\in\Z^d,\;  j\geq 0,  \; \ell\in\kauno:=\{1,\ldots, 2^d-1\},\\
\phi_{k}(x)&= \phi(x-\vk), \qquad \vk\in \Z^d,
\end{align*}
we can write any function $f\in L^2(\R^d)$ using the wavelet representation
\begin{align*}  
f(x)&=\sum_{\vk\in\Z^d}\langle f, \phi_\vk\rangle\phi_\vk(x)
+\sum_{\substack{j\geq 0,\ \vk\in \Z^d\\ \ell\in\kauno}}\langle f, \psi^\ell_{l j\vk}\rangle\psi^\ell_{ j\vk}(x)\\
&=\sum_{\vk\in\Z^d}f_{-1k}\phi_k(x)+\sum_{\substack{j\geq 0,\ \vk\in \Z^d
\\ \ell\in\kauno}}f^\ell_{ j\vk}\psi^\ell_{ j\vk}(x).
\end{align*}
Especially,  Daubechies wavelets with vanishing moments up to order $N-1\in\N$ are suitable for $r$-regular multi-resolution analysis if $r<0.1936(N-1)$. 
 Daubechies wavelets $\phi, \Psi^\ell$ ($\ell=1,\ldots, 2^d-1$)  of order $N$ are supported in the cube $[-L, L]^d$, where we have $L:=2N-1$. For more details about wavelets see \cite{Daubechies2006, Meyer1995}.


 We have the following equivalent representation to the Besov norm \eqref{eq:BesovNorm1} in case $p=q$
\begin{align}\label{equivalence}
\|f\|_{B^s_{pp}(\R^d)}
\simeq \left(\sum_{{j=-1}}^\infty 2^{jp\big(s+\frac d2-\frac dp\big)}
\|\textbf{f}_j\|_{\ell^p}^p\right)^{1/p},
\end{align}
where
\begin{align*}
\textbf{f}_j&=(f^\ell_{ j\vk})_{\vk\in \Z^d, \; \ell\in\kauno}\quad \text{with}\quad
\|\textbf{f}_j\|_{\ell^p}=\left(\sum_{\vk\in \Z^d, \; \ell\in\kauno}
|f^\ell_{ j\vk}|^p\right)^{1/p}  \qquad \textrm {if}\;\, j\geq 0,\quad \textrm{and}\\
\textbf{f}_{-1}&=(f_{-1\vk})_{\vk\in \Z^d} \quad \text{with}\quad
\|\textbf{f}_{-1}\|_{\ell^p}=\left(\sum_{\vk\in \Z^d}|f_{-1\vk}|^p\right)^{1/p}.
\end{align*}

\bigskip


\bigskip

We may reenumerate the full set of indices of our wavelets  
$$
\{(-1,k)_{k\in\Z^d}\}\cup \{ (\ell,j,k)_{\ell\in\kauno,\,j\in\N,\,k\in\Z^d}\}
$$ 
as $\{(\ell(l), j(l),k(l)):\ l\in \N\}$ (if $j(l)=-1$ we may set set $k(l)=-1$ and interpret $(-1,-1,k)=(-1,k)$). In the case $p=q$ we have the following trivial but essential consequence of \eqref{equivalence} which states that considering random Besov-space valued variables is abstractly  equivalent to considering random variables in a weighted $\ell^p$-space (or via another trivial isomorphism, in a non-weighted $\ell^p$-space).

\begin{proposition} Besov space $B_{pp}^s(\R^d)$  is isomorphic to the weighted $\ell^p$ space 
\begin{align*}
\ell^p_v=\{(a_l)_{l=1}^\infty:\ \|(a_l)\|^p_{\ell^p_v}=
\sum_{l=1}^\infty |v_la_l|^p<\infty\}
\end{align*}
where $v=(v_1,v_2,\dots),$ 
\begin{align}\label{v-weight} 
v_l=2^{j(l)\big(s+\frac d2-\frac dp\big)}.
\end{align}
The isomorphism is given by
\begin{align*}
(a_l)\mapsto \sum_{l=1}^\infty a_l\psi^\ell(l)_{j(l),\vk(l)}.
\end{align*}
\end{proposition}



For simplicity we will mainly consider the behaviour of random  functions on the unit cube $D=[0,1]^d$. Hence, while defining random functions via the wavelet decompositions we may welll set
$f^\ell_{ jk}=0$ if $\vk=(k_1,\dots,k_d)$
and  $|k_l|>L2^j-1$ for some $l=1,\dots,d$. Because of this,
we call the set of wavelet indices relevant to us an entire tree, which is given by
\begin{align*}  
\oT=\{(j,\vk)\in \N\times \Z^d\,,\, j\geq 0,\ \vk=(k_1,\dots,k_d),\ 
|k_l|\leq L2^j-1\}
\end{align*}
Note that we consider the same tree for all $\ell\in\kauno$.

We emphasize that in our theoretical results on the behaviour of the Besov prior we aim for modelling the generic local behaviour of the prior, which  in our setup  takes place \emph{in the interior of the unit cube $D=[0,1]^d$}. One may  of course fine-tune the definition of the prior suitably in connection with different boundary conditions, and study its behaviour also near the boundary, but since this is generally application specific we do not consider it in this paper.

Next we introduce random variables in Besov spaces, or equivalently in $\ell^p_v$. The generated measures are similar to the p-exponential measures whos consistency has been studied in \cite{Agapiou2018}. 

\begin{definition}\label{Def:BesovVariables}   
Let $\{\psi^\ell_{ jk}\}_{\ell\in\kauno,\, j\in\N,\,k\in\Z^d}$ be, as before, an $r$-regular wavelet basis for $L^2(\R^d)$. Let $(X^\ell_{j\vk})_{(j,\vk) \in\oT}$ be an i.i.d. sequence of real random variables with probability density function $\pi(x)\propto\exp (-\frac{1}{2\kappa^p}|x|^{p})$, in which case we say that $X^\ell_{jk}\sim \Nn_p(0,\kappa^p)$. 
Let the random function $f$ be defined as
\begin{align*}
f(x)=\sum_{\substack{(j,\vk)\in T\\ \ell\in\kauno}}h_j X^\ell_{ j,\vk}\psi_{j,\vk}(x),\quad x\in D,
\end{align*}
where $h_j=2^{-j\big(s+\frac d2-\frac dp\big)}$, $s<r$, are deterministic constants.  In light of \eqref{equivalence} say that $f$ is a $B^{s}_{pp}$-random variable.
\end{definition}

We note that random variables defined in Definition \ref{Def:BesovVariables} take values in Besov spaces $B^{t}_{pp}$, with $t<s-\frac{d}{p}$, a.s. and a realisation $f$ takes values in the space $B^{s-d/p}_{pp}$ only with probability zero. This follows directly from Theorem \ref{Thm:singsupp} below with $\gamma=d$. The space $B^{s}_{pp}$ plays similar role to the Cameron-Martin space  and informally $f$ has density proportional to $\exp(-\frac{1}{2\kappa^p}\|f\|_{B^s_{pp}}^p)$, see \cite{Agapiou2018, Lassas2009, Dashti2012}. If $p=2$ we obtain a Gaussian measure with Cameron-Martin space $B^s_{22}=H^s$. If $p=1$ $f$ is a semi-Laplace random variable, which means that $f$ is a hierarchical random variable which is determined in the last stage as a Laplace random variable.

\section{Generalisation to random functions having
singularities on random fractals}\label{Sec:Fractals}

In applications the strength of wavelets often appears in the fact that if large portion of small wavelet coefficients of a given function are replaced by zero, the new function corresponding to these truncated wavelet coefficients approximates well the original function. In particular, the singularities of the  original function are often preserved. Because of this we consider a model where a large part of wavelet coefficients are assumed to be zero. This idea has been previously studied e.g. by \cite{Abramovich1998, Chang2000, Jansen2012, Sendur2002} but here we emphasise the singularities by choosing the non-zero wavelet coefficients in a systematic way. 

We will consider the following set of trees in $\oT$, 
\begin{align*}  
\Gamma=\{T\subset\oT\,|\, 
\hbox{if $(j,\vk)\in T$ and $j\geq 1$ then $(j-1,[\vk/2])\in T$}\}
\end{align*} 
where $[\vk/2]=([k_1/2],\dots,[k_d/2])$ is the vector which
elements are  the integer parts of $k_l/2$. We call the set $\Gamma$ proper subtrees. 
The above definition means that if some node is in the tree $\Gamma$ then all of its parent nodes are also in the tree. That is, all the branches are connected to the root node.

\begin{definition}\label{Def:TreeRandomVariable} 
Let $\{\psi^\ell_{jk}\}_{\ell\in\kauno, j\in\N,\,k\in\Z^d}$ be an $r$-regular wavelet basis for $L^2(\R^d)$, and $\beta=2^{\gamma-d}$ with some $\gamma\in(-\infty,d\,]$. Consider pairs $(X^\ell,T)$ where $X^\ell$ is an $\R^{\oT}$-valued random variable, and $T\in {\Tt}$ is a random tree. We assume that $X$ and $T$ are independent random variables, having the following distributions
\begin{itemize}
\item The sequence $X^\ell$ consists i.i.d $X^\ell_{jk}\sim\Nn_p(0,\kappa^p)$, with probability density proportional to $\exp (-\frac{1}{2\kappa^p}|x|^{p})$, $\kappa>0$ and $1\leq p <\infty$.
\item The tree $T$ is build recursively
by choosing at each level new nodes into the tree with probability
$\beta$. More precisely, $T$ is determined with the following rule:
We assume that the root node $(0,0)$ is always chosen. For the rest of the levels;
Let $t_{j,\vk}$ be independent uniformly distributed random variables on $[0,1]$. When for a given $j$ all pairs  $(j,\vk)$ that are in $T$ are chosen, we choose pair $(j+1,l)$ to be in the tree $T$ if and only if $\big(j,[\frac l 2]\big)\in T$ and $t_{j+1,l}<\beta$.
\end{itemize}
Let $f$ be the random function 
\begin{align*}
f(x)=\sum_{\substack{(j,\vk)\in T\\ \ell\in\kauno}}h_{j} X^\ell_{j,\vk}\psi^\ell_{j,\vk}(x),\quad x\in D,
\end{align*}
where $h_{j}=2^{-j\big(s+\frac d2-\frac dp\big)}$, $s<r$, are deterministic constants. We say that $f$ is a $B^{s}_{pp}$-random variable with wavelet density $\beta$.
\end{definition}

Note that for a single $(\ell,j,k)$ the random variable $X^\ell_{j,\vk}\sim\mu=\Nn_p(0,\kappa^p)$ has probability distribution
\begin{align*}
X^\ell_{jk}\sim P_j\mu+(1-P_j)\delta_0,
\end{align*}
i.e. it vanishes with a probability $1-P_j$ (compare this with 
\cite{Abramovich1998}). 

Next we will employ theory of random fractals based on\cite{Dryakhlov2001}. Let us consider dyadic cubes
\begin{align*}
Q_{j,\vk}=\{y\in [0,1]^d:\ [jy]=\vk\}.
\end{align*}
For the random variable $T$ we define a random
fractal (a set-valued
random variable) 
\begin{align*}
\tau=\bigcap_{j=1}^\infty\bigcup_{\vk\in R_j}
\overline Q_{j,\vk}
\end{align*}
where $R_j=\{\vk:(j,\vk)\in T\}$ is set of level $j$ elements
in the random tree $T$.  The following theorem gives the Hausdorff dimension of the random fractal $\tau$, which is showed to coincide with the $r$ -singular support of the random variable $f$ in Theorem \ref{Thm:singsupp} below.

\begin{theorem}\label{Thm:DimOfTau}
Let $\beta =2^{\gamma-d}$, with $\gamma\in(-\infty,d\,]$, and $T$ be chosen as in Definition \ref{Def:TreeRandomVariable}. If $\gamma\in(-\infty,0]$ then $\tau$ is an empty set with probability one. If $\gamma\in(0,d]$ the set $\tau$ has Hausdorff dimension
\begin{align*}
\hbox{dim}_H\,(\tau)=\gamma,
\end{align*}
with probability $1-P_\beta$ and is empty with probability $P_\beta$, where $P_\beta$ is the solution to $P_\beta=((1-\beta)+\beta P_\beta)^{2^d}$. 
\end{theorem}
  
\begin{proof}
This result is known from the basic theory of Galton-Watson trees, see e.g. \cite[Section 8]{Karlin1975}, but for the readers benefit we sketch part of the argument here in the one-dimensional case.  We start by noting that $\tau$ is empty if and only if the tree $T$ terminates at some finite level, and we denote the probability for this by $P_\beta$. We will first look at the case $d=1$. Since the nodes of the tree are chosen recursively and a new node can only be chosen if its parent node is chosen we can write 
\begin{align*}
P_\beta=(1-\beta)^2+2(1-\beta)\beta P_\beta+\beta^2P_\beta^2=((1-\beta)+\beta P_\beta)^2. 
\end{align*}
Solving the above for $P_\beta$ gives us that the probability for the tree $T$ being finite, and for $\tau$ being an empty set, is $\left(\frac{1-\beta}{\beta}\right)^2$. 
For a general $d$ the result follows by a similar argument. We note that if $\beta\leq 2^{-d}$ there is no solution $P_\beta<1$ and the tree terminates almost surely at some finite level. 

Consider the random process $v(j)$ which is
the number of elements on $R_j=\{\vk:(j,\vk)\in T\}$.
Then $v(0)=1$ and the random variables $v(j)$, $j\geq 1$, follow binomial
distributions $v(j)\sim \hbox{Bin}(2^d v(j-1),\beta)$.
This  means that $v(j)$ is the number of offspring at level $j$ 
in the induced Galton-Watson branching process and when $\beta>2^{-d}$ we note that $\expec v(j) =(2^d\beta)^j>1$ for all $j\geq1$. 

We denote by $q_{jk}$ the ratio of the diameter of $Q_{jk}$ to the diameter of its parent when the parent is non-empty set, so that $q_{jk}=1/2$ with probability $\beta$ and $q_{jk}=0$ with probability $1-\beta$. We then notice that 
\begin{align*}
\E\bigg(\sum_{k=\{0,1\}^d}q_{1k}^\gamma\bigg)=2^d\left(\frac{1}{2}\right)^\gamma\beta=1
\end{align*}
which implies that when $\tau$ is non-empty its Hausdorf dimension is $\gamma$ a.s. \cite[Theorem 1.1]{Mauldin1986}.
\end{proof}

Consider next a realisation  of the random function $f$ corresponding to a realisation of $(X,T)=(X^\ell,T)_{\ell\in\kauno}$. Since the wavelets $\psi^\ell_{j,\vk}$ are in a H\"older space $C^r(\R^d)$ and the wavelet-representation is locally finite outside the  closed set $\tau$, we see that
\begin{align*}  
f|_{D\setminus \tau}\in C^r_{loc}(D\setminus \tau),\quad 
\hbox{i.e.}\quad \hbox{singsupp}_r(f)\subset \tau, 
\end{align*}
where $ \hbox{singsupp}_r(f)$ is the $C^r$-singular support of $f$. Motivated by this, we show next that $f$ is a Besov-space $B^{t}_{pp}$, $t<s-\gamma/p$, valued function and that the $C^r$-singular support
of $f$ is a.s. the random fractal $\tau$ corresponding to the realisation of the tree $T$.

\begin{theorem}\label{Thm:singsupp}
Let $f$ be a $B^{s}_{pp}$-random variable with wavelet density $\beta=2^{\gamma-d}$, with $\gamma\in(0,d]$, as in Definition \ref{Def:TreeRandomVariable}. Then, for all $t<s-\frac{\gamma}{p}$,
$f$ takes values in $B_{pp}^{t}$ almost surely  and 
$f\not\in B_{pp}^{s-\gamma/p}$ on the event
$$
G:=\{ \omega\in \Omega \; :\; \tau\not= \emptyset \}
$$
Moreover,  
\begin{align*}
\hbox{singsupp}_r(f)=\tau
\end{align*}
almost surely.
\end{theorem}

\begin{proof} 
We consider again the random process $v(j)\sim \hbox{Bin}(2^d v(j-1),\beta)$ which is the number of offspring at level $j$  in the induced Galton-Watson branching process.  Denote by $\bar{v}= \expec v(1) =2^d\beta>1$, and  $w(j)=v(j)/\bar{v}^j$. Next we consider the Besov norm of $f$, which according to \eqref{equivalence} has the same distribution as
\begin{equation}\label{eq:normi}
\| f\|_{B_{pp}^{t}}^p=\sum_{j=0}^{\infty}v(j)
2^{-jp(s-t)}\frac{\sum_{u=1}^{(2^d-1)v(j)}|Y_{j}^u|^p}{v(j)},
\end{equation}
where $Y_{j}^u$ are independent draws from $\Nn_p(0,\kappa^p)$.  

Let us first consider the case $t<s-\frac{\gamma}{p}$. Denote $\E(|Y_{j}^u|^p) = c_p$. Then 
\begin{align*}
\expec \bigg(\frac{\sum_{u=1}^{(2^d-1)v(j)}|Y_{j}^u|^p}{v(j)}\bigg)=(2^d-1)c_p
\quad\text{and}\quad
 \sum_{j=0}^\infty \bar{v}^j2^{-jp(s-t)}<\infty
\end{align*}
and hence we see
$$
\expec \|f\|^p_{B^{t}_{pp}} <\infty.
$$
This verifies that, for $t<s-\frac{\gamma}{p}$, we have $f\in B^{t}_{pp}$ almost surely.

Next we will consider the case $t=s-\frac{\gamma}{p}$.  
Since $\expec (v(1))^2 <\infty$, the basic theory of Galton-Watson
processes (see e.g. \cite{Athreya2004}) yields that the sequence 
$(w(j))_{j\geq 1}$ is an $L^2$-bounded martingale that has a limit $w_\infty$,
\begin{align*}
w_\infty=\lim_{j\to \infty}w(j)
\end{align*}
satisfying 
\begin{equation}\label{eq:rajamartingaali}
w_\infty \not=0\quad {\rm a.s.\;\; on}\;\; G.
\end{equation}
Naturally, $w_\infty =0$ in $ \Omega\setminus G$.
Since $w(j)$ is then a uniformly integrable martingale we 
have, by Doob's theorem, $L^1(\Omega)$-convergence
\begin{align}\label{L1-conv}
\lim_{j\to \infty}\|w(j)-w_\infty\|_{L^1(\Omega)}=0.
\end{align}

Write $Z_j:=\sum_{u=1}^{v(j)(2^d-1)}|Y_j^u|^p$.   
Let $m_j\in\{1,\dots,2^{jd}\}$ and  consider the random variable $\tilde Z_j$ that is the variable $Z_j$ conditioned on the set $\{\omega\, |\, v(j)=m_j\}$. Denoting $\text{Var}(|Y_{j}^u|^p) = c'_p$ we get $\expec (\tilde Z_j)=c_pm_j$ and $\text{Var}(\tilde Z_j) = c'_pm_j$. Hence we can conclude 
\begin{align*}
\expec (\tilde Z_j^2)=c'_pm_j+(c_pm_j)^2\leq c''_pm_j^2.
\end{align*}
Write $A:= \{\omega:\ \tilde Z_j\geq c_pm_j/2\}$ and let $a=\prob(A)$.
Then the Cauchy--Schwarz inequality yields
\begin{align*}
c_pm_j  =\expec \tilde Z_j 
&\leq (1-a)c_pm_j/2+\int_A\tilde Z_j\,d\prob\\
& \leq (1-a)c_pm_j/2+\Big(\int_\Omega\tilde Z_j^2 d\prob\Big)^{1/2}\Big(\int_\Omega\chi_A d\prob\Big)^{1/2}\\
& \leq (1-a)c_pm_j/2+m_j\sqrt{c''_pa}.
\end{align*}
Thus $c_p(1+a)\leq 2\sqrt{c''_pa}$ and we obtain that $a\geq a_0,$ where $a_0\in(0,1]$ depends only on $p$.
Hence for all $m_j\in\{1,\dots,2^{jd}\}$ it holds that
\begin{equation}\label{eq:alaraja}
\prob \big(\{Z_j\geq m_jc_p/2\}\; |\; v(j)=m_j\big)\geq a_0.
\end{equation}
One may observe that we essentially reproved a lemma due to Paley-Zygmund.

Let $q$ and $m_1,m_2,\ldots ,m_q$ be given positive integers with $m_i\in\{1,\dots,2^{id}\}$. This time we condition $Z_j$ on the set
$$
A_{m_1,m_2,\ldots ,m_q} =\{\omega\in \Omega \,:\, v(1)=m_1,\ldots ,v(q)=m_q\}.
$$
Notice that since the $Y_{j}^u$ are independent random variables, the estimate (\ref{eq:alaraja}) holds for $j\leq q$ as well if we replace there
the set $\{v(j)=m_j\}$ by $A_{m_1,m_2,\ldots ,m_q}.$, i.e.,
\begin{equation}\label{eq:alaraja 2}
\prob (\{Y_j\geq m_jc_p/2\}\; |\; A_{m_1,m_2,\ldots ,m_q})\geq a_0,\quad 
j\leq q.
\end{equation}
To apply this, consider variables $W_j:=Z_j/v(j)$. Observe that the variables $Z_j$, $1\leq j\leq q$, conditioned on the set $A_{m_1,m_2,\ldots ,m_q}$,  are independent. Thus (\ref{eq:alaraja 2}) implies that
$$
\prob \bigg(\Big\{\sum_{j=1}^q W_j <\frac {c_pq}2\Big\}\;\ \bigg |\; A_{m_1,m_2,\ldots ,m_q}
\bigg)
\leq (1-a_0)^{q}.
$$
Let $C>0$ and $\varepsilon>0$ be chosen arbitrarily. Choosing $q$ so large that $q>2Cc_p^{-1}$ and $(1-a_0)^{q}<\e$ we see that 
$$
\prob \bigg(\Big\{\sum_{j=1}^q W_j <M\Big\}\; |\; A_{m_1,m_2,\ldots ,m_q} \bigg)<\varepsilon.
$$
This does not depend on $q$ nor $m_1,m_2,\ldots ,m_q$ and we deduce that
$$  
\prob \bigg(\Big\{\sum_{j=1}^\infty W_j =\infty\Big\}\, \bigg|\,G\bigg)=1
$$
which implies,  in view of \eqref{eq:rajamartingaali} and the the expression \eqref{eq:normi} for the Besov norm, that almost surely $\prob (\{\|f\|^p_{B^{s-\frac{\gamma}{p}}_{pp}}=\infty\})=1$.

It remains to prove that ${\rm sing\, supp_r\,}(f) = \tau$ almost surely. By construction it is almost immediate that ${\rm sing\, supp_r\,}(f) \subset \tau$. Towards the other direction, given any dyadic subcube $Q$ of $D$, denote by $2Q$ the cube parallel to and with the same center as $Q$ but double the size of $Q$. Let $\varphi_Q$ denote a smooth cut-off function that is zero outside $2Q$ and one in a neighbourhood of $Q$. Because of the stochastic structure of the tree, the part of the tree (and the coefficients of our random Besov function) corresponding to $Q$ is (essentially) similar to the whole tree. Hence the above proof applies and we deduce that
$$
\textrm {a.s.} \quad  \|\varphi_Qf\|_{B_{pp}^{s-\gamma/p}}=\infty\quad \textrm{on the set}\quad \tau\cap \overline{Q}\not=\emptyset. 
$$
On the other hand, we have the equality of the following events:
$$
\{{\rm sing\, supp_r\,}(f) \subsetneq \tau\} =\bigcup_Q  \{\|\varphi_Qf\|_{B_{pp}^{r}}<\infty\}\cap \{\tau\cap\overline{Q}\not=\emptyset\},
$$
and the claim follows combining  these observations and the fact that the number of dyadic subcubes $Q$ is countable.
\end{proof}

%
%

\begin{corollary}
Let $f$ be a $B^s_{pp}$-random variable with wavelet density 
$\beta=2^{\gamma-d}$, $\gamma\in(-\infty,d]$, as in Definition \ref{Def:TreeRandomVariable}. Theorems \ref{Thm:DimOfTau} and \ref{Thm:singsupp} imply that for $\gamma\in(-\infty,0]$ a.s. $f\in C^r$, and for $\gamma\in(0,d]$ $\textup{singsupp}_r(f)$ is an empty set with probability $P_\beta$, with  $P_\beta$ being the solution of $P_\beta=((1-\beta)+\beta P_\beta)^{2^d}$, and $\dim_H(\textup{singsupp}_r(f))=\gamma$ with probability $1-P_\beta$. 
\end{corollary} 
  
\section{MAP-estimate for denoising problem}\label{Sec:exaples}  

\subsection{Discrete wavelet decomposition model}
In this section we study  signal and image denoising examples and show how the MAP estimator can be calculated explicitly. We start by introducing some notation for finite trees. For clarity we introduce the notation for the 1-dimensional case. The techniques readily generalise to 2-dimensional setting and the image denoising example is introduced at the end of the section.

We define an entire finite tree as a set 
\begin{equation*}   
\oT=\{(j,k)\in \N^2\ |\ 0\le j\le j_{\max},\ 0\le k\le 2^{j}-1\},
\end{equation*}
where $j_{\max}$ is a chosen terminating depth. We also consider finite proper subtrees
\begin{align*}  
\Gamma=\{T\in\oT\,|\, 
\hbox{if $(j,\vk)\in T$ and $j\geq 1$ then $(j-1,[\vk/2])\in T$}\},
\end{align*} 
where $[k/2]$ is the integer part of $k/2$.
We denote $(j',k')\trianglelefteq (j,k)$ if $(j,k)=(j',k')$ or $(j',k')\vartriangleleft(j,k)$ by which we mean that $(j,k)$ is an ancestor of $(j',k')$. We can then define the full subtree with a root node $(j,k)$ as
\begin{align*} 
\oT_{(j,k)} =\{(j',k')\in \N^2\ |\ (j',k')\trianglelefteq (j,k)\}.
\end{align*}
The size of a sub tree $\oT_{(j,k)}$ (number of the nodes) is ${s}_j=2^{j_{\max}-j+1}-1$. The parent of a node $(j,k)$ is $P(j,k)=(j-1,[k/2])$, and its left and right child are $C_0(j,k)=(j+1, 2k)$ and  $C_1(j,k)=(j+1, 2k+1)$ respectively.  If the nodes $(j,k)$ and $(j',k')$ have the same parent we say that they are siblings and denote $(j',k')\sim (j,k)$. 

We consider the denoising problem of recovering $f$ from
\begin{align*}
M=f+\Noise,        
\end{align*}  
where $\Noise=\sum w_{jk}\psi_{jk}$, $w_{jk}\sim \Nn(0,1)$, is a white noise process independent of $f$. We employ a discrete version of the random tree Besov prior introduced in the previous section and assume that, with some $T\in\Gamma$, $f$ can be written as
\begin{align}\label{eq:Unknown}  
\begin{split}  
f(x) & =\sum_{(j,k)\in T} \langle f,\psi_{jk}\rangle\psi_{jk}(x)\\
& =\sum_{(j,k)\in \oT} \tt_{jk} g_{jk}\psi_{jk}(x),
\end{split}  
\end{align} 
where $g_{jk} \sim \Nn(0,1)$ or $g_{jk} \sim \text{Laplace}(0,a)$, and $\tt_{jk}\in\{0,1\}$ defines if a node $(j,k)\in \oT$ is chosen i.e. term $g_{jk}\psi_{jk}(x)$ is non-zero. We assume that an independent node $t_{jk}$ is chosen with probability $\beta$, that is, $\Pp(t_{jk}=1)=\beta$ and $\Pp(t_{jk}=0)=1-\beta$. Then the sub tree $\Gamma$ contains the nodes for which $\tt_{jk}=1$, where $\tt_{jk}$ is defined as
\begin{align}\label{defb}
\tt_{jk}=\prod_{(j',k')\trianglerighteq(j,k)}t_{j'k'}.
\end{align}
Notice that this means that a coefficient can only be chosen if all of its ancestors have been chosen. Otherwise the coefficient is zero. All variables $g_{jk}, t_{jk}$ are assumed to be mutually independent. 

The data can be written as $M=\sum m_{jk}\psi_{jk}(x)$, where
\begin{align}\label{eq:Measurement}
m_{jk}=\left\{\begin{array}{ll}g_{jk}+w_{jk}, & \textrm{when}\ \tt_{jk}=1\\
w_{jk}, & \textrm{when}\ \tt_{jk}=0. \end{array} \right.
\end{align}  
The posterior distribution of $g=(g_{jk})_{(j,k)\in\oT}$ and $t=(t_{jk})_{(j,k)\in\oT}$, given data $m=(m_{jk})_{(j,k)\in\oT}$ can then be written as
\begin{align*}
\pi(g,t\ |\ m) & \propto\pi(m\ |\ g,t)\pi(g)\pi(t)\\
& =\prod_{(j,k)\in \oT} \pi(m_{jk}\ |\ g_{jk},\tt_{jk})\pi(g_{jk})\pi(t_{jk}).
\end{align*}

\subsection{Pruning and tree enforced soft thresholding algorithms}    
In this section we will show how the MAP estimator for the above denoising problem can be calculated explicitly. If we assume that  $g_{jk}\sim\Nn(0,1)$ the result is a pruning algorithm where the wavelet density $\beta<1/2$ acts as a regularisation parameter. Regularisation is achieved through turning branches of the wavelet tree on or off depending on  whether they are important for the reconstruction. If we assume $g_{jk}\sim\text{Laplace}(0,a)$ the outcome is a mixture of the above mentioned pruning algorithm and soft thresholding where the threshold is given by $a$.

We start  by noting that $g_{jk}$ has density 
\begin{align*}
\pi(g_{jk})\propto\exp(-R(g_{jk})),
\end{align*}
where $R(g_{jk})=\frac{g_{jk}^2}{2}$ when $g_{jk}\sim\Nn(0,1)$ and $R(g_{jk})=\frac{|g_{jk}|}{a}$ when $g_{jk}\sim\text{Laplace}(0,1)$.
We define
\begin{align*}
z_{jk}=\pi(m_{jk}\ |\ g_{jk},\tt_{jk})\pi(g_{jk})\pi(t_{jk})
\end{align*}
and note that since the root node $(0,0)$ is always chosen 
\begin{align*}
z_{0\,0}=\pi(m_{0\, 0}\ |\ g_{0\, 0},t_{0\, 0}=1)\pi(g_{0\, 0}).
\end{align*}
For general $(j,k)\in \oT$ the value of $z_{jk}$ depends on whether the node is chosen or not and we denote
\begin{align*}
z_{jk}^1 & =\pi(m_{jk}\ |\ g_{jk},\tt_{jk}=1)\pi(g_{jk})\pi(t_{jk}=1)
=\exp\Big(-\frac{1}{2}(m_{jk}-g_{jk})^2-R(g_{jk})\Big)\beta\\
z_{jk}^0 & =\pi(m_{jk}\ |\ g_{jk},\tt_{jk}=0)\pi(g_{jk})\pi(t_{jk}=0)
=\exp\Big(-\frac{1}{2}m_{jk}^2-R(g_{jk}\Big)(1-\beta).
\end{align*}

The problem of maximising $\pi(g,t\ |\ m)$ is equivalent to minimising  $\Big(-\log\big(\prod_{(j,k)\in \oT} z_{jk}\big)\Big)$ and we can write 
\begin{align*}
\min_{g,\,t}\bigg(-\log\Big(\prod_{(j,k)\in \oT} z_{jk}\Big)\bigg)
& =\min_{g,\,t}\bigg(-\log(z_{0\,0}) -\sum_{(j,k)\in \oT\setminus(0,0)}\Big(\tt_{jk}\log(z_{jk}^1)+(1-\tt_{jk})\log(z_{jk}^0)\Big)\bigg)\\
& =F(m).
\end{align*}
Notice that if $t_{jk}=0$ then $\tt_{j'k'}=0$ for all $(j',k')\trianglelefteq(j,k)$. This means that the whole branch is turned off and we attain minimum 
\begin{align*}
\min_{g}\bigg(-\sum_{(j',k')\in\oT_{(j,k)}}\log(z_{j'k'}^0)\bigg)=
\frac{1}{2}\|m|_{\oT_{(j,k)}}\|_{2}^2-{s}_j\log(1-\beta)
\end{align*}
by choosing $g|_{\oT_{(j,k)}}=0$. 

Next we denote the restriction of the function $F(m)$ to a sub tree $\oT_{(j,k)}$ by $F_{(j,k)}(m) = F(m|_{\oT_{(j,k)}})$.  We can then rewrite the minimisation problem in the following recursive form
\begin{align*}
F(m) & =\min_{g,\,t}\bigg(-\log(z_{0\,0})
-\sum_{k=0}^1 \Big(t_{1\,k}\big(\log(z_{1\,k}^1)
+\sum_{(j',k')\in \oT_{(1,k)}\setminus(1,k)} \log(z_{j'k'})\big)\\ 
&\quad\,+(1-t_{1\,k})\sum_{(j',k')\in\oT_{(1,k)}}\log(z_{j'k'}^0)\Big)\bigg)\\
& =\min_{t_{1\,0},t_{1\,1}}\bigg(-\log(z_{0\,0})+\sum_{k=0}^1 \Big(t_{1\,k}\big(F_{(1,k)}(m)-\log\beta\big)\\
&\quad\, +(1-t_{1\,k})\Big(\frac{1}{2}\|m|_{\oT_{(1,k)}}\|^2-{s}_1\log(1-\beta)\Big)\Big)\bigg).
\end{align*}
We can repeat the above recursion until we reach the lowest level $j_{max}$. 

To solve the minimisation problem we start by calculating the  minimum for $-\log\big(\pi(m_{jk}\ |\ g_{jk},\tt_{jk}=1)\pi(g_{jk}))$, $j=j_{\max}$. If $g_{jk}\sim\Nn(0,1)$ then 
\begin{align*}
\min_{g_{jk}}\Big( -\log\big(\pi(m_{jk}\ |\ g_{jk},\tt_{jk}=1)\pi(g_{jk})\big)\Big)
=\min_{g_{jk}} \Big(\frac{1}{2}(m_{jk}-g_{jk})^2+\frac{1}{2}g_{jk}^2\Big)
=\frac{1}{4}m_{jk}^2,
\end{align*}
and the minimum is attained by choosing $g_{jk}=m_{jk}/2$. If $g_{jk}\sim\text{Laplace}(0,a)$ then 
\begin{align*}
\min_{g_{jk}}\Big( -\log\big(\pi(m_{jk}\ |\ g_{jk},\tt_{jk}=1)\pi(g_{jk})\big)\Big)
& = \min_{g_{jk}} \Big(\frac{1}{2}(m_{jk}-g_{jk})^2+\frac{1}{a}|g_{jk}|\Big)\\
& =\begin{cases}
      \frac{1}{2}m_{jk}^2 & \text{when}\ |m_{jk}|\leq a \\
     \frac{1}{a}\big(m_{jk}-\frac{1}{2a}\big) & \text{when}\ m_{jk}> a\\
     -\frac{1}{a}\big(m_{jk}+\frac{1}{2a}\big) & \text{when}\ m_{jk}<- a
    \end{cases}
\end{align*}
and
\begin{align*}
\argmin_{g_{jk}} \Big(\frac{1}{2}(m_{jk}-g_{jk})^2+\frac{1}{a}|g_{jk}|\Big)
=\begin{cases}
      0, & \text{when}\ |m_{jk}|\leq a \\
     m_{jk}-a & \text{when}\ m_{jk}> a\\
     m_{jk}+a & \text{when}\ m_{jk}<- a.
 \end{cases}
\end{align*}

After finding the minimising values for the nodes on the bottom level $j_{\max}$ we move one level up and test if choosing a bottom node, while paying the penalty $-\log\beta$ for doing so, gives a smaller value than not choosing it. The same procedure is then carried out on every level. See Algorithm \ref{Al:MAP} below for case $g_{jk}\sim\Nn(0,1)$. 

\newpage
\begin{algorithm}[H]\label{Al:MAP}
{$j=j_{\max}$}\\
\For{$k = 0$ \KwTo $2^{j}-1$}{  
\[ F_{(j,k)}(m)=\frac{1}{4}|m_{jk}|^2;\]}

\While {$1\leq j\leq j_{\max}-1$}{
$j=j-1$\\
\For{$k = 0$ \KwTo $2^{j}-1$}{
\begin{align*}
 F_{(j,k)}(m) =
\frac{1}{4}|m_{jk}|^2
& + \sum_{i=0}^1\min_{\substack{t_{C_i(j,k)}} }
\bigg(t_{C_i(j,k)}\Big(F_{C_i(j,k)}(m)-\log \beta\Big)\\
& +(1-t_{C_i(j,k)})\Big(\frac{1}{2} \|m|_{\oT_{C_i(j,k)}}\|^2-\overline{s}_{j+1}\log(1-\beta)\Big) \bigg);
\end{align*}
}}
Let $t$ be the minimising sequence used to achieve $F(m)=F_{(0,0)}(m)$.\\
\uIf{$\tilde{t}_{jk}=1$}{
\begin{align*}
g_{jk}=m_{jk};
\end{align*}}
\Else{
\begin{align*}
g_{jk}=0;
\end{align*}}
\caption{Pseudocode for finding the minimising $t$ and $g$ recursively.}
\end{algorithm}
\vspace{3mm}

Notice that when $g_{jk}\sim\Nn(0,1)$ the denoising problem is regularised only by turning of branches that do not carry enough information. If $g_{jk}\sim \text{Laplace}(0,a)$ regularisation is achieved through soft thresholding, where the threshold depends on the prior variance, and by excluding branches with small wavelet coefficients. 

\subsection{Signal and image denoising examples}

\begin{figure}[H]
\centering
\begin{subfigure}{.5\textwidth}
  \centering
  \includegraphics[trim={1.5cm 0cm 0.5cm 0.5cm },clip,height=5cm, width=8cm]{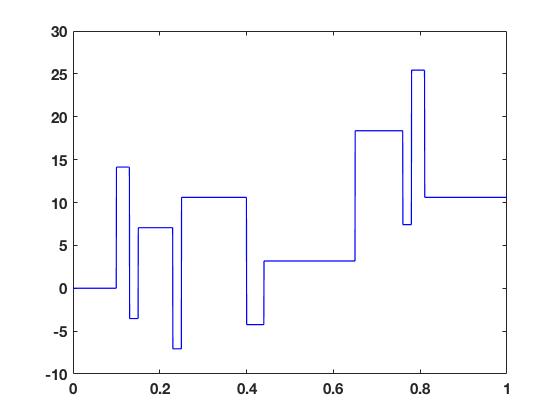}
\end{subfigure}%
\begin{subfigure}{.5\textwidth}
  \centering
  \includegraphics[trim={1.5cm 0cm 0.5cm 0.5cm },clip, height=5cm, width=8cm]{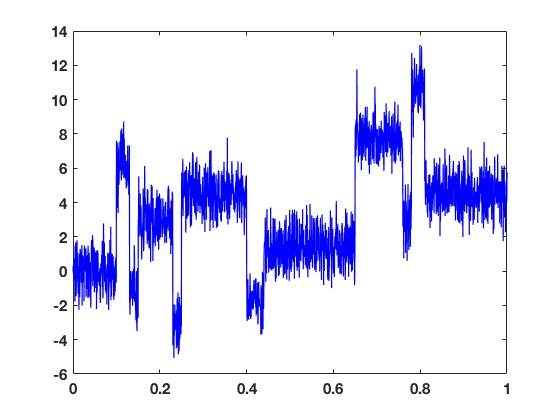}
\end{subfigure}
\caption{Original blocks signal on left  and the observed noisy signal, with signal to noise ratio 3, on right.}
\label{Fig:BlocksData}
\end{figure}

\begin{figure}[H]
\centering
\begin{subfigure}{.5\textwidth}
  \centering
  \includegraphics[trim={1.5cm 0cm 0.5cm 0.5cm },clip,height=5cm, width=8cm]{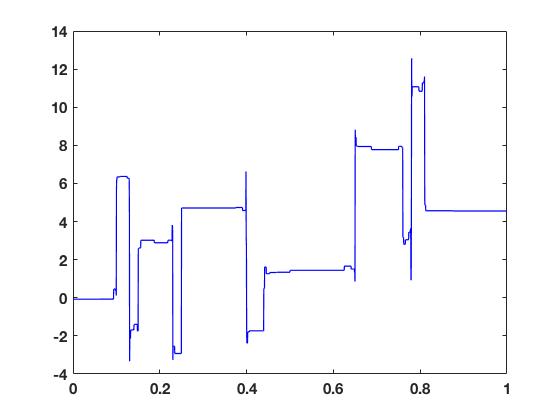}
   \includegraphics[trim={2.5cm 0cm 0.5cm 0.5cm },clip,height=4cm, width=6cm]{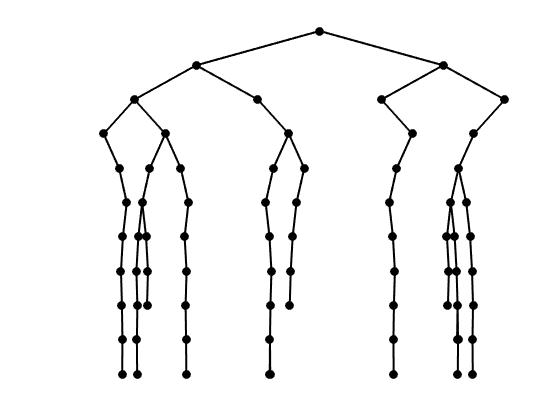}
\end{subfigure}%
\begin{subfigure}{.5\textwidth}
  \centering
  \includegraphics[trim={1.5cm 0cm 0.5cm 0.5cm },clip, height=5cm, width=8cm]{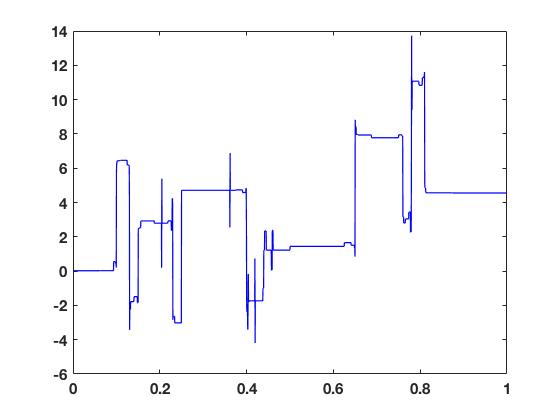}
  \includegraphics[trim={2.5cm 0cm 0.5cm 0.5cm },clip, height=4cm, width=6cm]{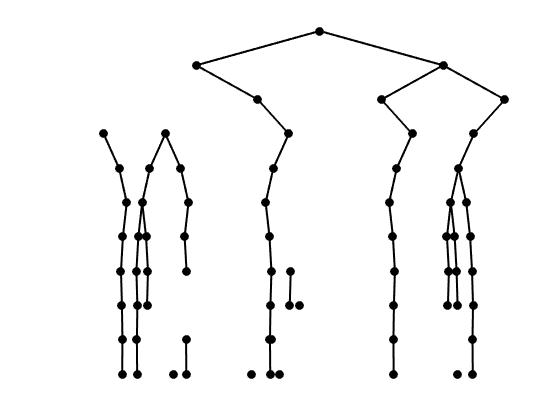}
\end{subfigure}
\caption{The denoised signal using pruning algorithm on top left and the denoised signal that was attained using hard thresholding on top right. Below them are the wavelet trees corresponding to the estimators.}
\label{Fig:DenoisedBlocks}
\end{figure}

In our first example we consider the blocks test data from \cite{Donoho1994} which is displayed in Figure \ref{Fig:BlocksData}. We assume that only a noisy signal, with noise ratio 3, is observed and want to denoise it. We employ semi-Gaussian random tree Besov prior with Haar wavelets, in which case the MAP estimator is given by the pruning algorithm.  We also tested denoising the signal with several Matlab denoising packages of which hard thresholding performed the best. Figure \ref{Fig:DenoisedBlocks} shows that the pruning algorithm performs better than hard thresholding. The $\ell^2$ error and root mean squared error  between the original and our denoised signal are $9.5$ and $0.21$ respectively. For the thresholded signal the errors are $13.7$ and $0.30$. The denoised signals and the corresponding wavelet trees are shown in Figure \ref{Fig:DenoisedBlocks} and one can clearly see how choosing the proper tree model is beneficial for the reconstruction.

\begin{figure}[H]
\centering
\begin{subfigure}{.33\textwidth}
  \centering
  \includegraphics[trim={2cm 0.7cm 0cm 0cm },clip,height=3.5cm, width=5.8cm]{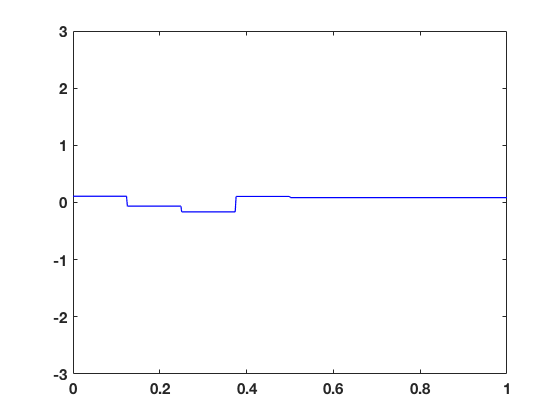}
  \includegraphics[trim={2cm 0.7cm 0cm 0cm },clip, height=3.5cm, width=5.8cm]{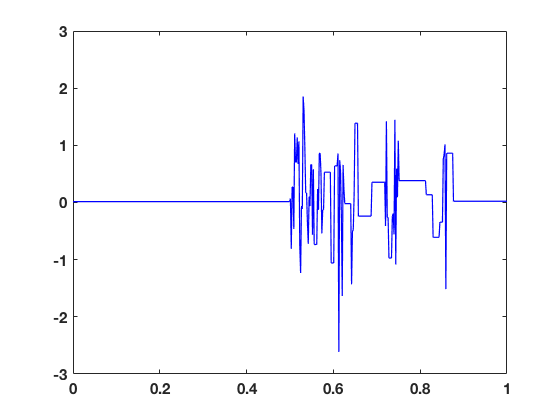}
  \includegraphics[trim={2cm 0.7cm 0cm 0cm },clip, height=3.5cm, width=5.8cm]{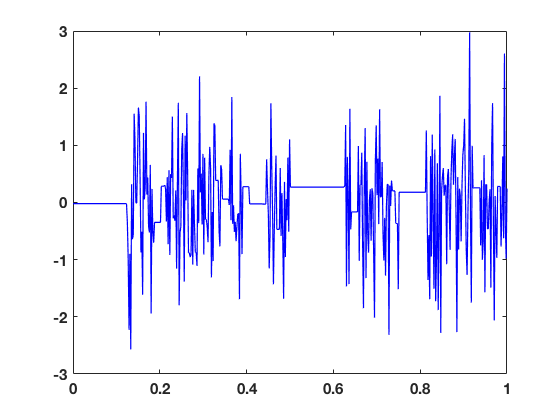}
\end{subfigure}%
\begin{subfigure}{.33\textwidth}
  \centering
  \includegraphics[trim={2cm 0.7cm 0cm 0cm },clip, height=3.5cm, width=5.8cm]{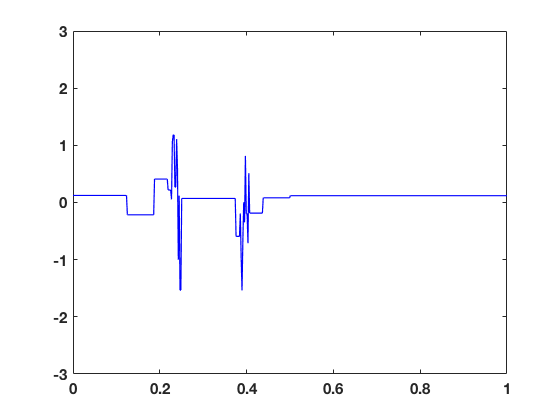}
  \includegraphics[trim={2cm 0.7cm 0cm 0cm },clip, height=3.5cm, width=5.8cm]{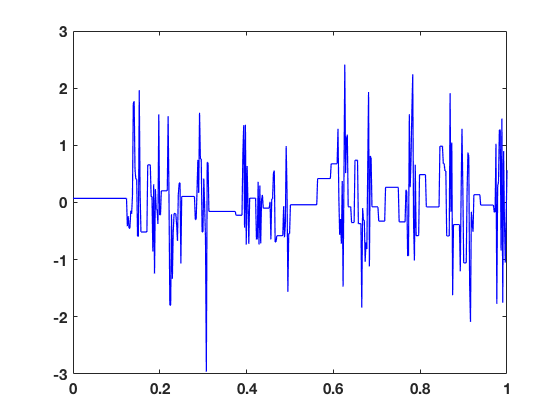}
  \includegraphics[trim={2cm 0.7cm 0cm 0cm },clip, height=3.5cm, width=5.8cm]{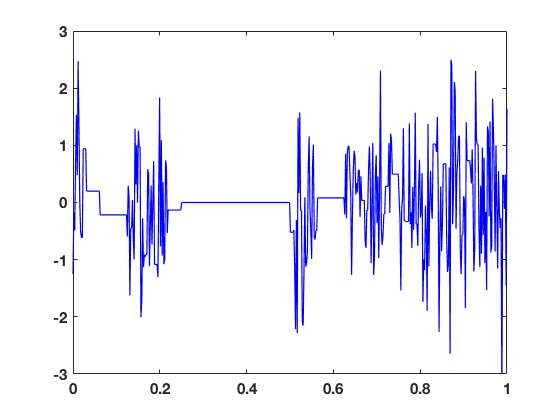}
\end{subfigure}
\begin{subfigure}{.33\textwidth}
  \centering
  \includegraphics[trim={2cm 0.7cm 0cm 0cm },clip, height=3.5cm, width=5.8cm]{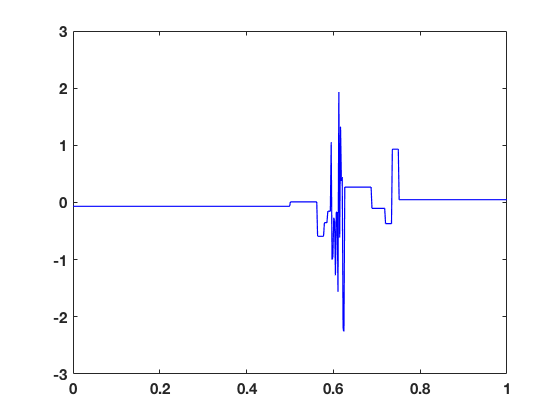}
  \includegraphics[trim={2cm 0.7cm 0cm 0cm },clip, height=3.5cm, width=5.8cm]{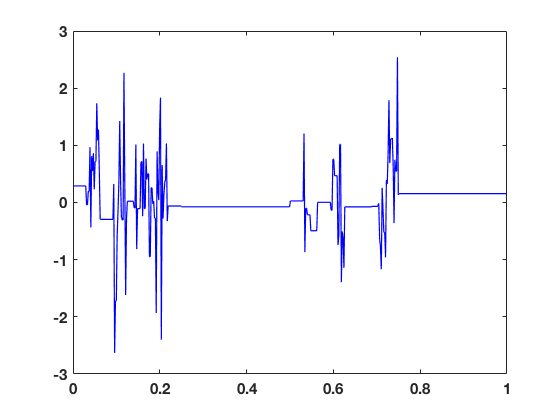}
  \includegraphics[trim={2cm 0.7cm 0cm 0cm },clip, height=3.5cm, width=5.8cm]{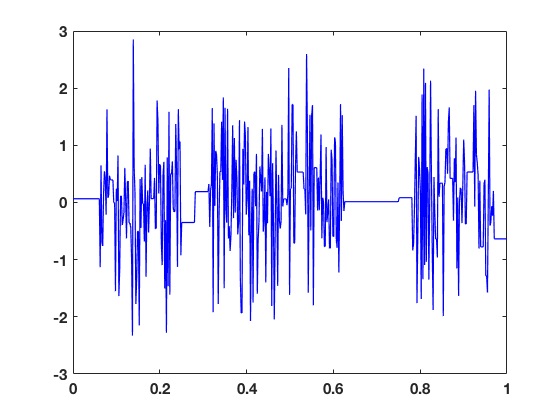}
\end{subfigure}
\caption{Prior draws from the 1-dimensional semi-Gaussian random tree Besov prior with Haar wavelets. The wavelet densities are $\beta=0.6$ on the first row, $\beta=0.75$ on the middle row and $\beta=0.9$ on the bottom row.}
\label{Fig:1dPriorDraws}
\end{figure}

In Figure \ref{Fig:1dPriorDraws} we present some prior draws  from semi-Gaussian random tree Besov priors with Haar wavelets with different wavelet densities $\beta$. As expected, when $1/2<\beta<1$ is small the realisations are mostly flat with some small areas of large jumps, while with large $\beta$ the realisations are more erratic. 

\begin{figure}[H]
\centering
\begin{subfigure}{.5\textwidth}
  \centering
  \includegraphics[trim={1.5cm 0cm 0.5cm 0.5cm },clip,height=6cm, width=8cm]{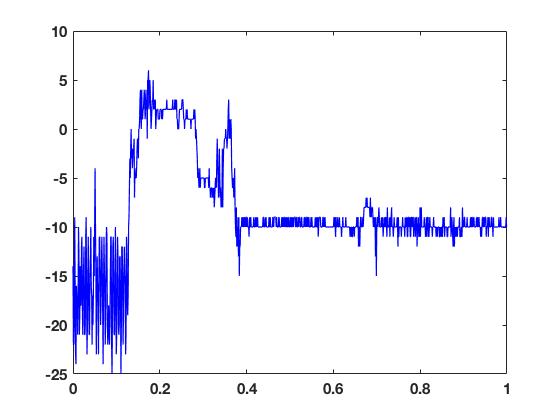}
\end{subfigure}%
\begin{subfigure}{.5\textwidth}
  \centering
  \includegraphics[trim={1.5cm 0cm 0.5cm 0.5cm },clip, height=6cm, width=8cm]{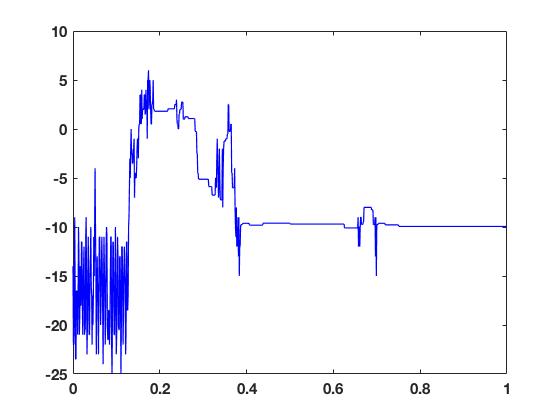}
\end{subfigure}
\caption{Measured accelerometer data on left and the denoised signal on right. }
\label{Fig:1dData}
\end{figure}

Next we consider real accelerometer data collected from a human subject with a wearable device. The accelerometer erroneously records slight movement even when the device is still and our aim is to use the pruning algorithm with Haar wavelets to denoise the signal. The original and denoised signal are presented in Figure \ref{Fig:1dData}. The pruning algorithm returns a denoised signal where the areas of large acceleration have been left untouched while the parts where the device was still have been cleaned and are flat as they should be. 

\begin{figure}[H]
\centering
\begin{subfigure}{.5\textwidth}
  \centering
  \includegraphics[trim={2cm 0cm 3cm 1cm },clip,height=7.5cm]{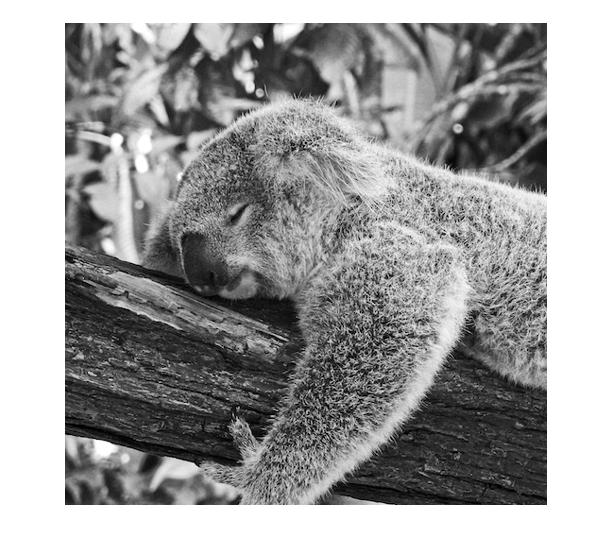}
\end{subfigure}%
\begin{subfigure}{.5\textwidth}
  \centering
  \includegraphics[trim={2cm 0cm 3cm 1cm },clip, height=7.5cm]{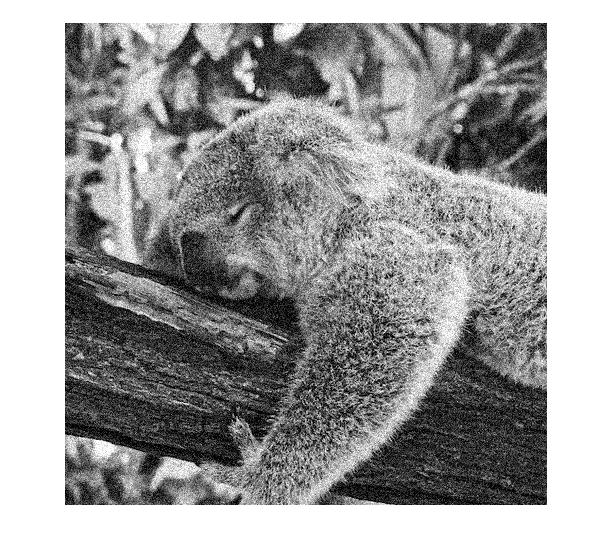}
\end{subfigure}
\caption{Original sharp image and image with Gaussian noise with variance $0.01$.}
\label{Fig:2dData}
\end{figure}

\begin{figure}[H]
\centering
\begin{subfigure}{.33\textwidth}
  \centering
  \includegraphics[trim={2cm 0cm 2cm 0cm },clip, width=1\linewidth]{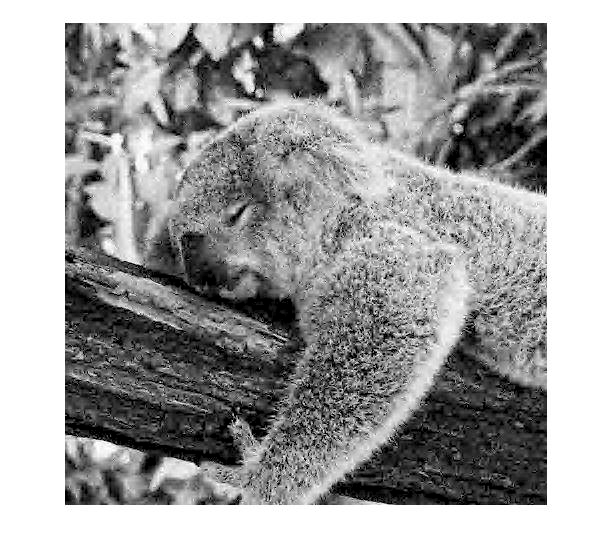}
\end{subfigure}%
\begin{subfigure}{.33\textwidth}
  \centering
  \includegraphics[trim={2cm 0cm 2cm 0cm },clip, width=1\linewidth]{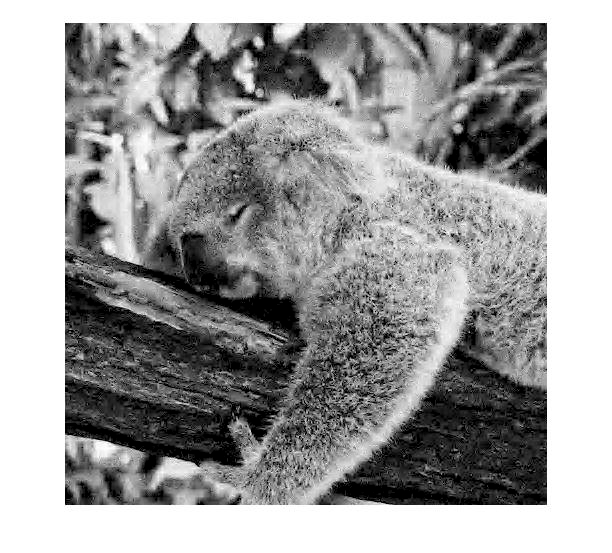}
\end{subfigure}
\begin{subfigure}{.33\textwidth}
  \centering
  \includegraphics[trim={2cm 0cm 2cm 0cm },clip, width=1\linewidth]{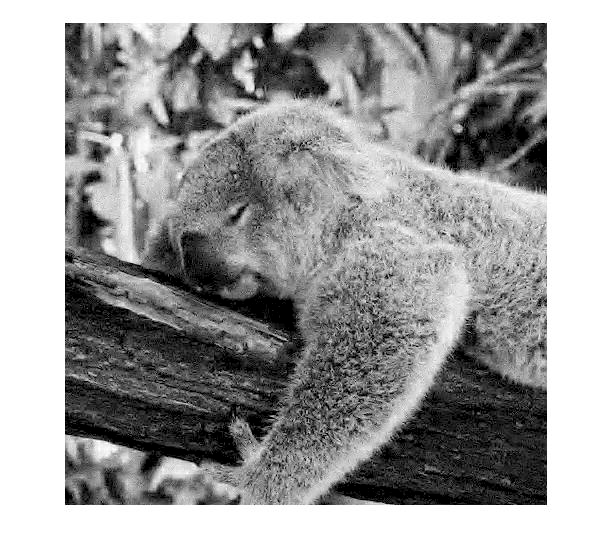}
\end{subfigure} 
\caption{From left to right; Denoised images attained using pruning, tree enforced soft thresholding and soft thresholding algorithms. The peak signal-to-noise ratios are $23.0$, $23.7$ and $23.6$, and the structural similarity indices are $0.69$, $0.71$ and $0.69$ respectively.}
\label{Fig:KoalaDenoised}
\end{figure}

In our final example we study image denoising. The original and noisy image are shown in Figure \ref{Fig:2dData}. We consider the MAP estimators arising from semi-Gaussian and semi-Laplace random tree Besov priors with Daubechies 2 wavelets, which are given by pruning and tree enforced soft thresholding algorithms respectively. The image has been extended by mirroring to avoid boundary effects, and cut down to the original size after using the algorithm. As expected the tree enforced soft thresholding performs better than the pruning algorithm since it allows two type of regularisation;  through pruning and soft thresholding. The peak signal-to-noise ratio (using the original image as a reference) for the pruned image is $23.0$, and $23.7$ for the image given by the tree enforced soft thresholding. The structural similarity indices are $0.69$ and $0.71$ respectively. We also denoised the image using soft thresholding with Daubechies 2 wavelets and the peak signal-to-noise ratio for the best achieved reconstruction is $23.6$ and structural similarity index is $0.69$. The three reconstructions are presented in Figure \ref{Fig:KoalaDenoised}.

\begin{figure}[H]
\centering
\begin{subfigure}{.33\textwidth}
  \centering
  \includegraphics[trim={4cm 2cm 2cm 0cm },clip, width=1\linewidth]{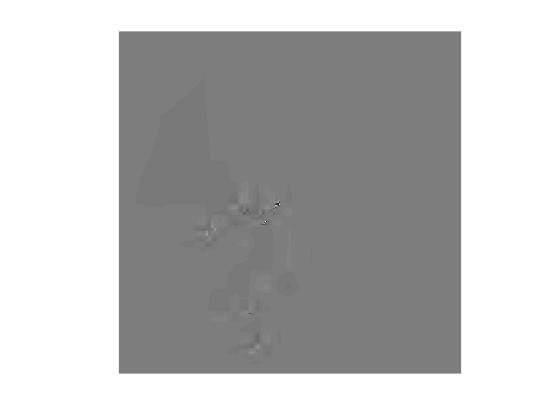}
  \includegraphics[trim={4cm 2cm 2cm 0cm },clip, width=1\linewidth]{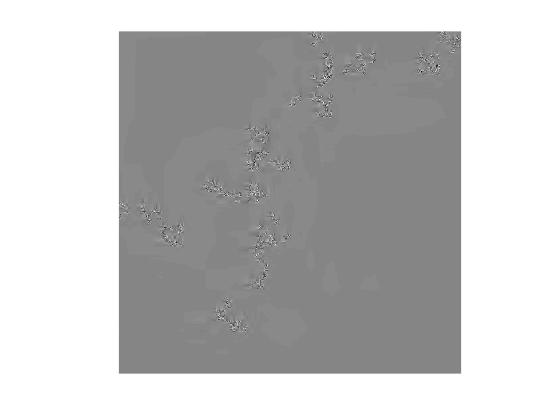}
  \includegraphics[trim={4cm 2cm 2cm 0cm },clip, width=1\linewidth]{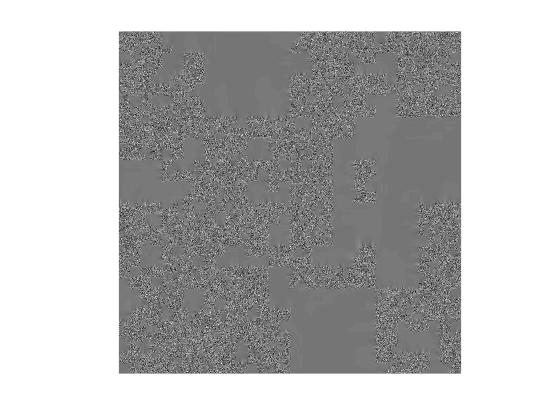}
\end{subfigure}%
\begin{subfigure}{.33\textwidth}
  \centering
  \includegraphics[trim={4cm 2cm 2cm 0cm },clip, width=1\linewidth]{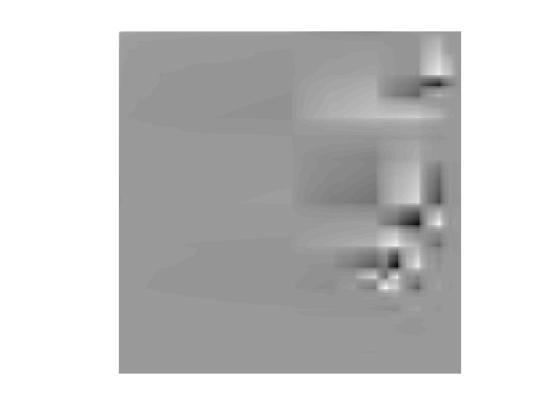}
  \includegraphics[trim={4cm 2cm 2cm 0cm },clip, width=1\linewidth]{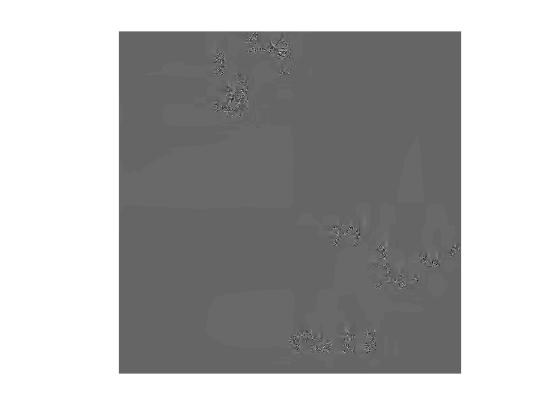}
  \includegraphics[trim={4cm 2cm 2cm 0cm },clip, width=1\linewidth]{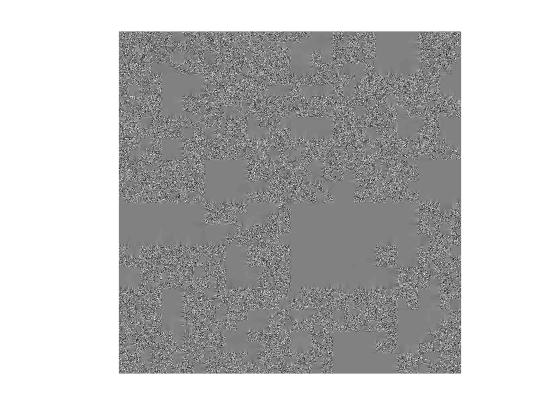}
\end{subfigure}
\begin{subfigure}{.33\textwidth}
  \centering
 \includegraphics[trim={4cm 2cm 2cm 0cm },clip, width=1\linewidth]{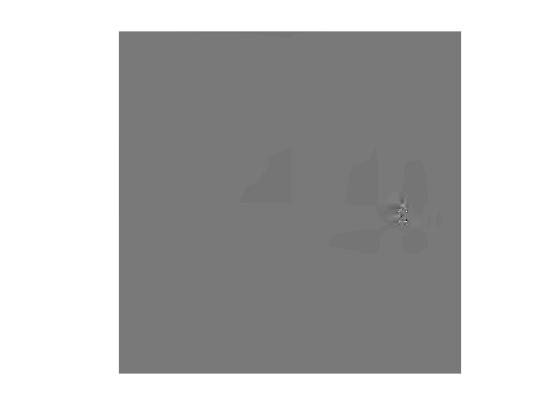}
  \includegraphics[trim={4cm 2cm 2cm 0cm },clip, width=1\linewidth]{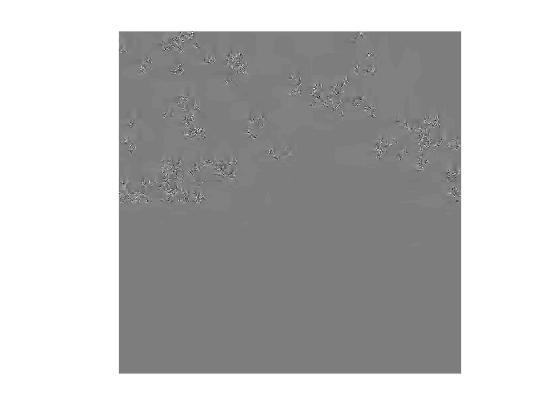}
  \includegraphics[trim={4cm 2cm 2cm 0cm },clip, width=1\linewidth]{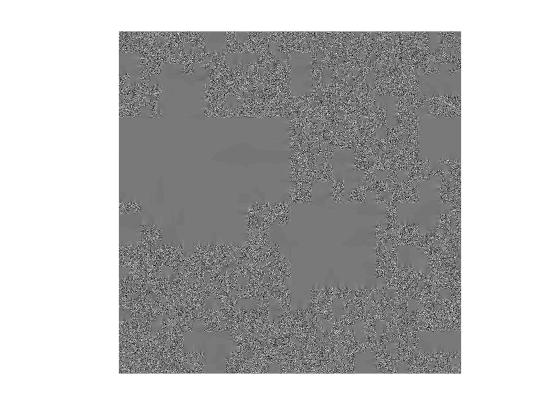}
\end{subfigure}
\caption{Prior draws from the 2-dimensional semi-Gaussian random tree Besov prior with Daubechies 2 wavelets. The wavelet densities are $\beta=0.3$ on the first row, $\beta=0.6$ on the middle row and $\beta=0.9$ on the bottom row.}
\label{Fig:2dPriorDraws}
\end{figure}

Random draws from the semi-Gaussian prior with Daubechies 2 wavelets and different wavelet densities can be found in Figure \ref{Fig:2dPriorDraws}.

\newpage
\subsubsection{Some results for general linear inverse problems}
We will next show that when a wavelet branch does not carry enough information it will be turned off also in the case of a general linear forward operator $A$.  We assume the random tree Besov prior introduced in \eqref{eq:Unknown}. We can then write
\begin{align*}
\pi(f\ |\ M) &\propto \pi(M-Af\g f)\pi(g)\pi(t)\\
& = \exp\Big(-\frac{1}{2}\|M-Af\|^2\Big)\exp\big(-R(g)\big)\beta^{n(t)}(1-\beta)^{(s_1-{n(t)})}
\end{align*}
where ${n(t)}$ is the number of elements in a tree $\Gamma=\{(j,k)\in\oT\ |\ \tt_{jk}=1\}$. We are interested in the minimisation problem 
\begin{align*}
\min_f \big(-\log\pi(f\ |\ M)\big)
=\min_{g,\,t}\Big(\frac{1}{2}\|M-Af\|^2+R(g)-n(t)\log\beta-(s_1-n(t))\log(1-\beta)\Big).
\end{align*}

We will next show that if the measurement $m|_{T_{(j,k)}}$ has small enough $\ell^2$ norm in some $T_{(j,k)}$ the subtree will be turned off.

\begin{lemma}\label{Lem:finite 2}
Let $M=Af+\Noise$, where $\Noise$ is a centred Gaussian noise process. We assume a random tree Besov prior defined in \eqref{eq:Unknown}. If $\|m|_{T_{(j,k)}}\|^2<\e$, where $0<\e<\min\{\log(1/\beta-1),1\}$, $0<\beta<1/2$ then the tree maximising the posterior $\pi(f\ |\ M)$ is a subset of $\oT\setminus\oT_{(j,k)}$.  
\end{lemma}

\begin{proof}
The result follows directly from the fact that if $\|m\|^2<\e$ the tree that maximises the posterior is an empty tree $T=\emptyset$. For an empty tree $n(t)=0$ and we have
\begin{align*}
\min_g \big(-\log\pi(f\ |\ M)\big)=\frac{1}{2}\|m\|^2-s_1\log(1-\beta).
\end{align*}
On the other hand, if $n(t)\geq1$ we can write
\begin{align*}
\min_{g,t} \big(-\log\pi(f\ |\ M)\big) 
& \geq \min_{t}\big(-n(t)\log\beta-(s_1-n(t))\log(1-\beta)\big)\\
& > \min_{t}\big(n(t)(\e-\log(1-\beta))-(s_1-n(t))\log(1-\beta)\big)\\
& >\|m\|^2-s_1\log(1-\beta)
\end{align*}
which concludes the proof. 
\end{proof}

%
%
%
%
%
%
%
%
%

\addcontentsline{toc}{chapter}{Bibliography}


\paragraph{Acknowledgements.}
H.K. acknowledges partial support by Cantab Capital Institute for the Mathematics of Information (RG83535). M.L. and S.S. were supported by the Academy of Finland through the Finnish Centre of Excellence in Inverse Modelling and Imaging 2018-2025, decision number 312339. E.S was supported by the Finnish Academy grant 1309940.

\bibliographystyle{siam}
\bibliography{Inverse_problems_references}

\end{document}